\newcommand{\siamyesno}[2]{#2}   
\siamyesno{
\documentclass[onetabnum,onefignum,nohypdvips,final]{siamart171218}
\usepackage{amssymb,amsmath,epsfig,verbatim,enumitem}
\newtheorem{remark}[theorem]{Remark}
\newtheorem{ass}[theorem]{Assumption}
}{
\documentclass[11pt,a4paper]{article}
\usepackage{amssymb,amsmath,amsthm,epsfig,verbatim,xcolor,enumitem}
\newtheorem{theorem}{\sc Theorem.}[section]
\newtheorem{lemma}[theorem]{\sc Lemma.}
\newtheorem{remark}[theorem]{\sc Remark.}

\renewcommand{\theequation}{\arabic{section}.\arabic{equation}}
\newenvironment{AMS}%
{{\upshape\bfseries AMS subject classifications. }\ignorespaces}{}
\newenvironment{keywords}{{\upshape\bfseries Key words. }\ignorespaces}{}

}
\usepackage[normalem]{ulem}
\usepackage{ifpdf}  
\ifpdf
  \DeclareGraphicsExtensions{.eps,.pdf,.png,.jpg}
\else
  \DeclareGraphicsExtensions{.eps,.ps}
\fi

\newcommand{\RZ}{{\bR} \slash {\mathbb Z}}
\newcommand{\bR}{{\mathbb R}}

\newcommand{\drho}{\;{\rm d}\rho}
\newcommand{\ds}{\;{\rm d}s}
\newcommand{\du}{\;{\rm d}u}
\newcommand{\dt}{\;{\rm d}t}
\newcommand{\Id}{I\!d}
\newcommand{\nabs}{\nabla_{\!s}}

\newcommand{\dd}[1]{\frac{\rm d}{{\rm d}#1}}
\newcommand{\ddt}{\dd{t}}
\newcommand{\ratio}{{\mathfrak r}}
\newcommand{\Vh}{\underline{V}^h}
\newcommand{\xpihx}{\sigma}
\newcommand{\yRhy}{\varsigma}
\newcommand{\pihxxh}{e^h}
\newcommand{\Rhyyh}{{\mathfrak e}^h}

\def\epsilon{\varepsilon} 
\def\hat{\widehat}

\makeatletter
\renewcommand{\uuline}{%
  \bgroup
  \UL@setULdepth
  \markoverwith{%
    \lower\ULdepth\hbox{%
      \kern-.03em%
      \vtop{%
        \hrule width.2em%
        \kern 0.6pt 
        \hrule
      }%
      \kern-.03em%
    }%
  }%
  \ULon
}
\makeatother
\siamyesno{}{
\textwidth 455pt \oddsidemargin 0pt \evensidemargin 0pt \headsep
0pt \headheight 0pt \textheight 655pt \parskip 10pt \parindent 0pt

}

\begin{document}
\title{
Finite element schemes with tangential motion
for fourth order geometric curve evolutions in arbitrary codimension}
\author{Klaus Deckelnick\footnotemark[2]\ \and 
        Robert N\"urnberg\footnotemark[3]}

\renewcommand{\thefootnote}{\fnsymbol{footnote}}
\footnotetext[2]{Institut f\"ur Analysis und Numerik,
Otto-von-Guericke-Universit\"at Magdeburg, 39106 Magdeburg, Germany \\
{\tt klaus.deckelnick@ovgu.de}}
\footnotetext[3]{Dipartimento di Mathematica, Universit\`a di Trento,
38123 Trento, Italy \\ {\tt robert.nurnberg@unitn.it}}

\date{}

\maketitle

\begin{abstract}
We introduce novel finite element schemes for curve diffusion and 
elastic flow in arbitrary codimension.
The schemes are based on a variational form of a system that includes a
specifically chosen tangential motion. We derive optimal
$L^2$-- and $H^1$--error bounds for 
continuous-in-time semidiscrete finite element approximations that use
piecewise linear elements.
In addition, we consider fully discrete schemes and, in the case of curve
diffusion, prove unconditional stability for it.
Finally, we present several numerical simulations, including some
convergence experiments that confirm the derived error bounds. The presented
simulations suggest that the tangential motion leads to equidistribution
in practice.
\end{abstract} 



\renewcommand{\thefootnote}{\arabic{footnote}}

\siamyesno{
\pagestyle{myheadings}
\thispagestyle{plain}
\markboth{K. DECKELNICK AND R. N\"URNBERG}
{FINITE ELEMENT SCHEMES WITH TANGENTIAL MOTION FOR FOURTH ORDER FLOWS}
}{}

\setcounter{equation}{0}
\section{Introduction} 

Curves endowed with a geometric energy, and gradient flows 
related to these energies, are of interest in differential
geometry, materials science, continuum mechanics, biological applications,
microelectronics and computer vision, 
\cite{Mullins57,LinS04,GoyalPL05,TuO-Y08,DorflerN19}. 
Two of the most frequently studied energies are the length functional
\begin{equation} \label{eq:L}
L(\Gamma) = \int_\Gamma 1 \ds
\end{equation}
and the elastic energy
\begin{equation} \label{eq:E}
E(\Gamma) = \tfrac12 \int_\Gamma |\varkappa|^2 \ds + \lambda L(\Gamma).
\end{equation}
Observe that \eqref{eq:E} 
combines a bending energy involving the modulus of the curvature vector
$\varkappa$ of the curve $\Gamma$ with a length contribution that 
penalizes growth if $\lambda>0$. In this paper we consider the case of closed
curves evolving in $\bR^d$ with $d\geq2$ arbitrary.
The most natural gradient flows for \eqref{eq:L} and \eqref{eq:E} are their
respective $L^2$--gradient flows, giving rise to curve shortening flow,
\cite{GageH86,Grayson87},
and elastic flow (also called curve straightening flow),
\cite{Wen95,Polden96,DziukKS02}.
In the case of \eqref{eq:L} the $H^{-1}$--gradient flow
is also of interest, and leads to curve diffusion,
\cite{TaylorC94,ElliottG97a,EscherMS98,GigaI98,GigaI99,DziukKS02,curves3d}. 
In the planar case, $d=2$, this flow is often also called surface diffusion 
for curves and has the important
property of conserving the area enclosed by the curve. In fact,
surface diffusion was proposed in \cite{Mullins57} as an
evolution law for a free surface enclosing a solid phase, which changes its
shape due to the diffusion of atoms along the surface. Later a derivation in
the context of rational thermodynamics was given in \cite{DaviG90}. \\
In this paper, we investigate the numerical approximation 
of curve diffusion and elastic flow. It turns out that
both of these fourth order geometric evolution equations are closely related,
and so our numerical analysis can be applied to both cases with only minor
modifications. From now on, and throughout this paper, we consider
a parametric description of the family of evolving curves
$(\Gamma(t))_{t\in[0,T]}$, $T>0$. In particular, we assume that
$\Gamma(t) = x(I,t)$, where $I = \RZ$ is the periodic unit interval and
$x : I\times [0,T] \to \bR^d$. On letting $\cdot_s$ denote differentiation 
with respect to arclength, the unit tangent along the curve is 
$\tau=x_s$. Then $x$ parameterizes a family of curves moving by curve 
diffusion if it satisfies the partial differential equation
\begin{equation} \label{eq:Vcd}
P x_t = - \nabs^2 \varkappa,
\end{equation}
where $P = \Id - \tau \otimes \tau$ and $\nabs \phi= P \phi_s$ for a vector 
field $\phi$. Similarly, elastic flow is described by solutions to
\begin{equation} \label{eq:Vel}
P x_t = - \nabs^2 \varkappa - \tfrac12 |\varkappa|^2 \varkappa
+ \lambda \varkappa.
\end{equation}
Observe that the tangential velocity, $x_t \cdot \tau$, is not prescribed in
\eqref{eq:Vcd} and \eqref{eq:Vel}, which is a natural consequence of the fact
that reparameterizations of the evolving family of curves does not affect
the geometric evolution.  \\
It is the aim of this paper to introduce a suitable
tangential motion in \eqref{eq:Vcd} and \eqref{eq:Vel} in such a way that
the obtained system of partial differential equations 
\begin{itemize}[noitemsep,topsep=0pt]
\item is strictly parabolic,
\item asymptotically exhibits solutions that are nearly arclength 
parameterizations,
\item admits a weak formulation that is amenable to discretization by finite
elements with a corresponding error analysis.
\end{itemize}
Our approach will be first developed for the curve diffusion flow. The main
idea is to require the tangential motion to be such that $x$ not only satisfies
\eqref{eq:Vcd}, but in addition reduces its Dirichlet energy 
$\int_I | x_\rho |^2 \drho$ in time, thus driving the parameterization towards 
one that is proportional to arclength. This idea can be viewed as a natural
extension of similar approaches for the curve shortening flow 
(\cite{DeckelnickD95,ElliottF17}) to the fourth order problem.
Intuitively it is clear that on the discrete level this will yield nice
distributions of mesh points, since large deviations in the length element
$|x_\rho|$ are penalized. In fact, in practice we do observe 
asymptotically equidistributed discretizations for the schemes introduced and
analysed in this paper. \\
Let us briefly review the existing literature on numerical
approximations of parametric formulations for curve diffusion and elastic flow,
with a particular emphasis on the available error analysis.
Various numerical schemes for curve diffusion have been proposed in
\cite{DziukKS02,BanschMN05,MikulaS05,triplej,curves3d,fdfi,BaoZ21,JiangL21,BaoL25}.
To the best of our knowledge, no error bounds for such schemes have been
proved so far.
But we mention the results in \cite{BanschMN04,DeckelnickDE03} for the error
analysis of surface diffusion, for the case of hypersurfaces evolving in
$\bR^3$, in the context of graph formulations.
With regards to elastic flow, numerical methods have been proposed in 
\cite{DziukKS02,triplej,willmore,DeckelnickD09,curves3d,pwf,Bartels13a,PozziS23}.
In \cite{DeckelnickD09} error estimates are shown, 
while \cite{Bartels13a} contains a partial convergence result 
for a scheme that approximates elastic flow of inextensible curves.
Here we recall that the numerical methods in \cite{DziukKS02,DeckelnickD09} 
discretize a purely normal flow, which can lead to very nonuniform meshes and
coalescence of vertices in practice. The BGN-schemes from the series of papers
\cite{triplej,willmore,curves3d,pwf,fdfi}, on the other hand, enjoy nearly
uniform distributions of mesh points in practice, see also 
\cite{BaoZ21,JiangL21} for related contributions and \cite{MikulaS05,DuanL24} for alternative
approaches involving tangential motion.
But for now there is no error analysis for these schemes. 
A first step to combine error analysis and good mesh properties for fourth
order flows was recently achieved in \cite{PozziS23}.
There a gradient flow for an energy functional like \eqref{eq:E}, but
with the length of the curve replaced by the Dirichlet energy
$\int_I | x_\rho |^2 \drho$, was introduced and analysed.
Since the Dirichlet energy is not invariant with respect to
repara\-meterization, the resulting flow is no longer geometric.
Nevertheless, this modified flow has the same stationary points as elastic
flow.
For more details on discretization methods for geometric partial differential
equations we refer to the two review articles \cite{DeckelnickDE05,bgnreview}. \\
We end this section with a few comments about notation. 
Throughout, $C$ denotes a generic positive constant independent of 
the mesh parameter $h$.
At times $\epsilon$ will play the role of a (small)
positive parameter, with $C_\epsilon>0$ depending on $\epsilon$, but
independent of $h$.

\setcounter{equation}{0}
\section{Mathematical formulation} \label{sec:mf}

Consider a family $(\Gamma(t))_{t\in[0,T]}$ of evolving curves that are given by $\Gamma(t) = x(I,t)$, where
$ x : I \times [0,T] \to \bR^d$ satisfies $| x_\rho| > 0$ in $I \times [0,T]$. Then the unit tangent on $\Gamma(t)$,
the curvature vector of $\Gamma(t)$ and the 
orthogonal projection onto the normal space of $\Gamma(t)$ are given by the 
following identities in $I$, see e.g.\ \cite{bgnreview}:
\begin{equation*} 
\tau = x_s = \frac{x_\rho}{| x_\rho |}, \quad
\varkappa =\frac{\tau_\rho}{| x_\rho |} =
\frac{1}{| x_\rho |}\bigl(\frac{x_\rho}{| x_\rho |} \bigr)_\rho, \quad 
P=\Id-\tau \otimes \tau.
\end{equation*}
With a view towards introducing variational discretizing methods for
curve diffusion, \eqref{eq:Vcd}, and elastic flow, \eqref{eq:Vel}, we follow 
the common approach of rewriting these fourth order problems in terms of two 
second order problems for the position vector $x$ and a further variable $y$. 
Rather than making the frequent choice $y=\varkappa$, in this paper we take 
\begin{equation} \label{eq:y}
y = \frac{x_{\rho\rho}}{|x_\rho|^2} ,
\end{equation}
so that 
\begin{equation*} 
\varkappa = \frac{x_{\rho \rho}}{| x_\rho |^2} - \frac{x_{\rho \rho} \cdot x_\rho}{| x_\rho |^3} \frac{x_\rho}{| x_\rho |} =  Py. 
\end{equation*}
The idea to use \eqref{eq:y} as a second variable extends the approach from
\cite{DeckelnickD95}, where the formulation 
$x_t = \frac{x_{\rho\rho}}{|x_\rho|^2}$ was used for curve 
shortening flow, to the fourth order flows considered in this paper. \\
Observing that  $P_s = -\tau_s \otimes \tau - \tau \otimes \tau_s = - P y \otimes \tau - \tau \otimes P y$, we calculate
\begin{equation*} 
\nabs \varkappa = (P y)_s - \bigl( (Py)_s \cdot \tau \bigr) \tau = P \bigl( P y_s +P_s y  \bigr) = P y_s - (\tau \cdot y) Py,
\end{equation*}
and hence 
\begin{align} \label{eq:kss}
\nabs^2 \varkappa& = P \bigl( P y_{ss} + P_s y_s  \bigr) - P \bigl( (\tau_s \cdot y) Py + (\tau \cdot y_s) Py + (\tau \cdot y) P_s y + (\tau \cdot y) Py_s \bigr) \nonumber \\
& = P y_{ss}  -2 ( \tau \cdot y_s) Py    - | Py |^2 Py +  (\tau \cdot y)^2 Py -( \tau \cdot y) P y_s \nonumber \\
&= P y_{ss}  -2 ( \tau \cdot y_s) Py    - | y |^2 Py + 2 (\tau \cdot y)^2 Py -( \tau \cdot y) P y_s .
\end{align}
Moreover, we have on noting \eqref{eq:y} that
\begin{equation} \label{eq:yss}
y_s = \frac{y_\rho}{| x_\rho |}, \quad y_{ss} = \frac1{|x_\rho|} \left( \frac{y_\rho}{|x_\rho|} \right)_\rho
= \frac1{|x_\rho|^2} y_{\rho\rho} - (y \cdot \tau) \frac{y_\rho}{|x_\rho|}.
\end{equation}
Combining \eqref{eq:kss} and \eqref{eq:yss} yields that
\begin{align} \label{eq:kssyss}
\nabs^2 \varkappa & = P \frac{y_{\rho\rho}}{|x_\rho|^2}
- 2(\frac{y_\rho}{|x_\rho|} \cdot \tau) Py
- 2 (y \cdot \tau)  P \frac{y_\rho}{|x_\rho|} 
+ 2(y \cdot \tau)^2 P y - |y|^2 Py \nonumber \\
& = \frac1{|x_\rho|^2} P \left[
y_{\rho\rho} -  2(y_\rho \cdot x_\rho) y 
- 2 (y \cdot x_\rho) y_\rho
+ 2(y \cdot x_\rho)^2 y - |x_\rho|^2|y|^2 y \right] .
\end{align}

\subsection{Curve diffusion} 
Our aim now is to introduce a system in such a way that \eqref{eq:Vcd} holds,
while the tangential part is chosen suitable to give desired properties. Starting from \eqref{eq:kssyss}, the most general ansatz
would hence be to consider the system
\begin{equation} \label{eq:xtansatz}
x_t = \frac1{|x_\rho|^2} \left( 
- y_{\rho\rho} +  2(y_\rho \cdot x_\rho) y 
+ 2 (y \cdot x_\rho) y_\rho
- 2(y \cdot x_\rho)^2 y + |x_\rho|^2|y|^2 y \right) 
+ \alpha \frac{x_\rho}{|x_\rho|^2} ,
\end{equation}
whose solutions clearly satisfy \eqref{eq:Vcd}, i.e.\ 
$P x_t  = -\nabs^2 \varkappa$ and therefore $\ddt \int_I |x_\rho| \drho \leq 0$ as a parameterization of the $H^{-1}$--gradient flow
of length.
Motivated by analogous properties for the DeTurck trick applied to curve
shortening flow, see \cite{DeckelnickD95}, we now attempt to find a coefficient $\alpha$ such that
solutions to \eqref{eq:xtansatz} satisfy in addition
\begin{equation} \label{eq:dtE}
\ddt \int_I |x_\rho|^2 \drho \leq 0.
\end{equation} 
To this end, we compute with the help of \eqref{eq:y} and \eqref{eq:xtansatz} 
that
\begin{align*} 
\tfrac12 \ddt \int_I |x_\rho|^2 \drho 
& = \int_I x_\rho \cdot x_{\rho,t} \drho 
= -\int_I y \cdot x_t |x_\rho|^2 \drho \nonumber \\ & 
= -\int_I |y_\rho|^2 \drho
- 2 \int_I y_\rho \cdot x_\rho |y|^2 \drho
- 2 \int_I (y \cdot x_\rho) y_\rho \cdot y \drho \nonumber \\ & \quad
+ 2 \int_I (y \cdot x_\rho)^2 |y|^2 \drho
- \int_I |x_\rho|^2|y|^4 \drho 
- \alpha \int_I x_\rho \cdot y \drho \nonumber \\ & 
= -\int_I |y_\rho|^2 + 2 y_\rho \cdot x_\rho |y|^2 + |x_\rho|^2|y|^4 \drho
\nonumber \\ & \quad
- \int_I y \cdot x_\rho [ 2 y_\rho \cdot y - 2 y \cdot x_\rho |y|^2
+ \alpha ] \drho \nonumber \\ & 
= -\int_I | y_\rho  + |y|^2 x_\rho|^2 \drho 
- \int_I y \cdot x_\rho [ 2 y_\rho \cdot y - 2 y \cdot x_\rho |y|^2
+ \alpha ] \drho.
\end{align*}
Choosing
\begin{equation*} 
\alpha = 2 y \cdot x_\rho |y|^2 - 2 y_\rho \cdot y
\end{equation*}
we find that solutions of the system
\begin{align} \label{eq:xt}
x_t & = \frac1{|x_\rho|^2} \left( 
- y_{\rho\rho} +  2(y_\rho \cdot x_\rho) y 
+ 2 (y \cdot x_\rho) y_\rho
- 2(y \cdot x_\rho)^2 y + |x_\rho|^2|y|^2 y \right) 
\nonumber \\ & \quad
+ \frac2{|x_\rho|^2} \left( y \cdot x_\rho |y|^2 - y \cdot y_\rho \right)
x_\rho
\end{align}
satisfy \eqref{eq:Vcd} and have the desired property \eqref{eq:dtE}. 
In particular, solutions to \eqref{eq:xt} satisfy the identity
\begin{equation} \label{eq:dtE0}
\ddt \int_I |x_\rho|^2 + \int_I | y_\rho  + |y|^2 x_\rho|^2 \drho  = 0.
\end{equation}
Let us rewrite some of the terms on the right hand side of \eqref{eq:xt} in a 
more convenient form, namely
\begin{align*} 
& 
(y \cdot x_\rho) y_\rho - (y \cdot x_\rho)^2 y + (y \cdot x_\rho) |y|^2 x_\rho
-(y \cdot y_\rho)x_\rho  \nonumber \\ & \quad
=  (y_\rho \otimes x_\rho) y - (x_\rho \otimes y_\rho) y +(x_\rho \cdot y) \bigl( (x_\rho \otimes y)y- (y \otimes x_\rho) y \bigr).
\end{align*}
Inserting this relation into \eqref{eq:xt} we obtain
\begin{equation} \label{eq:xt2}
|x_\rho|^2 x_t = - y_{\rho\rho} + F_{cd}(x_\rho,y,y_\rho) y,
\end{equation}
where $F_{cd}(a,b,c) \in \bR^{d \times d}$ is given by $F_{cd}=F_1+F_2$ with
\begin{subequations} \label{eq:defF}
\begin{align}
F_1(a,b,c) &= \bigl( 2 a \cdot c + | a |^2 | b |^2  \bigr)\Id, \label{eq:defF1} \\
F_2(a,b,c) &= 2 \bigl( c \otimes a -  a \otimes c \bigr) + 2 a \cdot b \bigl( a \otimes b - b \otimes a \bigr), \label{eq:defF2}
\end{align}
\end{subequations}
which corresponds to a splitting of $F_{cd}$ into a symmetric
and an anti-symmetric part. In particular, it holds that
\begin{equation} \label{eq:antisymm}
F_2(a,b,c) z \cdot z = 0 \qquad \mbox{ for all } a,b,c,z \in \bR^d.
\end{equation}
In the planar case, $d=2$, one has
\begin{displaymath}
F_2(a,b,c)b  = 2 \bigl( a^\perp \cdot c - (a \cdot b) a^\perp \cdot b \bigr) b^\perp,
\end{displaymath}
with $\cdot^\perp$ denoting the anti-clockwise rotation through $\frac{\pi}{2}$.

Let us briefly give formal reasons why we expect the above chosen tangential 
motion to have a positive effect on the discretization. Firstly,
as indicated in the introduction, the fact that the Dirichlet energy of $x$
satisfies the estimate \eqref{eq:dtE} means that $|x_\rho|$ should in general 
only show small variations. 
Moreover, the diffusive term in \eqref{eq:dtE0} 
gives some additional information.
To make this more precise, 
we temporarily assume that the solution $x$ of \eqref{eq:xt} exists 
globally in time and satisfies $0< | x_\rho | \leq C_0$ on $I \times [0,\infty)$. Then we deduce from \eqref{eq:dtE0} that
\begin{displaymath}
 \int_0^\infty \int_I | (y_\rho + |y|^2 x_\rho ) \cdot x_\rho |^2 \drho \dt \leq C_0^2
\int_0^\infty \int_I |y_\rho + |y|^2 x_\rho|^2 \drho \dt \leq \frac{C_0^2}{2} \int_I | x_\rho(\cdot,0) |^2 \drho < \infty,
\end{displaymath}
so that we expect $ (y_\rho + |y|^2 x_\rho) \cdot x_\rho$  to be  small for large $t$. Since
\begin{displaymath}
 (y_\rho + |y|^2 x_\rho) \cdot x_\rho  =   y_\rho \cdot x_\rho + y \cdot  | x_\rho |^2 y 
   =  y_\rho \cdot x_\rho + y \cdot x_{\rho \rho} 
 =   (y \cdot x_\rho)_\rho = \bigl( \frac{x_{\rho \rho} \cdot x_\rho}{| x_\rho |^2} \bigr)_{\rho}  = \bigl( \log | x_\rho | \bigr)_{\rho \rho} 
\end{displaymath}
this will in turn imply that 
$ \log | x_\rho|$ and hence also 
$ | x_\rho |$ itself will be nearly constant for large $t$ giving rise
to an almost arclength parameterization.
It is then natural to expect that discretizations based on \eqref{eq:xt} 
will lead to almost uniform distributions of grid points, and this 
is indeed what we will observe in the numerical experiments.

\subsection{Elastic flow}

On recalling that $\varkappa =Py$, we may write
\begin{displaymath}
- \tfrac12 | \varkappa |^2 \varkappa + \lambda \varkappa 
= - \tfrac12 | Py |^2 Py + \lambda Py 
=  \frac{1}{| x_\rho|^2} P \left[ - \tfrac12 \bigl( | x_\rho |^2 |y|^2 - (y \cdot x_\rho)^2 \bigr)y + \lambda | x_\rho |^2 y \right].
\end{displaymath}
If we combine this relation with  \eqref{eq:xt2}, we see that solutions to 
\begin{equation} \label{eq:xt2el}
|x_\rho|^2 x_t = - y_{\rho\rho}  +  F_{el}(x_\rho,y,y_\rho)y 
\end{equation}
satisfy \eqref{eq:Vel} if we choose $F_{el}(a,b,c)=F_{cd}(a,b,c) + F_3(a,b)$, 
where
\begin{equation} \label{eq:defF3}
F_3(a,b) = \bigl(
- \tfrac12 \bigl( |a|^2 |b|^2- (a \cdot b)^2 \bigr) + \lambda |a|^2
\bigr) \Id ,  \qquad  a,b \in \bR^d.
\end{equation}
Clearly, \eqref{eq:xt2el} differs from \eqref{eq:xt2} only in the additional
term $F_3(x_\rho,y)y$ on the right hand side. 
Note, however, that \eqref{eq:dtE0}, and the more general \eqref{eq:dtE},
will in general no longer hold. Consequently, the
heuristic argument regarding the asymptotic arclength parameterizations
shown at the end of the previous section no longer applies. 
Nevertheless we will see in
the numerical experiments that our approach works very well also in the case of elastic flow.

\setcounter{equation}{0}
\section{Discretization} \label{sec:sd}

We saw in the previous section that a solution of the system
\begin{subequations} \label{eq:sys}
\begin{align} 
| x_\rho |^2 x_t + y_{\rho \rho} & = F(x_\rho,y,y_\rho)y, \label{eq:sys1} \\
|x_\rho |^2 y - x_{\rho \rho} & = 0 \label{eq:sys2}
\end{align}
\end{subequations}
satisfies \eqref{eq:Vcd} if $F=F_{cd}$, and \eqref{eq:Vel} if $F=F_{el}$. 
Note that by inserting \eqref{eq:sys2} into \eqref{eq:sys1}, the system has the
form
\[
x_t = - \frac{x_{\rho\rho\rho\rho}}{|x_\rho|^4} + \text{lower order terms},
\]
so that we obtain a strictly parabolic problem. \\
In what follows we assume that \eqref{eq:sys}
has a unique solution $(x, y)$ with
$x \in C^1([0,T];H^4(I,\bR^d))$ and $x_{tt} \in L^2(0,T;L^2(I,\bR^d))$
satisfying 
\begin{equation*} 
| x_\rho | \geq c_0,  \quad 
| x_\rho | + | y|  \leq C_0 \quad \mbox{ on } I \times [0,T], \quad  \int_0^T | x_t |_{1,\infty}^2 \dt \leq C_0
\end{equation*} 
for constants $0<c_0<C_0$. Note that this implies 
$y\in C^1([0,T];H^2(I,\bR^d))$. \\[2mm]
In order to define a semidiscrete approximation, we choose
 a partition $0=q_0 < q_1< \ldots < q_{J-1} < q_J=1$ of $[0,1]$ into 
intervals $I_j=[q_{j-1},q_j]$ and set  $h_j= q_j-q_{j-1}$ as well as $h=\max_{j=1,\ldots,J} h_j$.
We assume that there exists a positive constant $c$ such that
\begin{equation*} 
h \leq c h_j, \quad 1 \leq j \leq J,
\end{equation*}
so that the resulting family of partitions of $[0,1]$ is quasi-uniform. 
Within $I$ we identify $q_J=1$ with $q_0=0$ and define the finite element 
space
\begin{displaymath}
V^h = \{\chi \in C^0(I) : \chi\!\mid_{I_j} \mbox{ is affine},\ 
j=1,\ldots, J\}, \quad \Vh = [V^h]^d.
\end{displaymath}
Let $\{\chi_j\}_{j=1}^J$ denote the standard basis of $V^h$.
For later use, we let $\pi^h:C^0(I)\to V^h$ 
be the standard Lagrange interpolation operator, and similarly
for $\pi^h:C^0(I,\bR^d)\to \Vh$.
It is well-known that for 
$k \in \{ 0,1 \}$, $\ell \in \{ 1,2 \}$ and $p \in [2,\infty]$ it holds that
\begin{subequations}
\begin{alignat}{2}
h^{\frac 1p - \frac 1r} \| \eta \|_{0,r} 
+ h | \eta |_{1,p} & \leq C \| \eta \|_{0,p} 
\qquad && \forall\ \eta \in V^h, \qquad r \in [p,\infty], 
\label{eq:inverse} \\
| z - \pi^h z |_{k,p} & \leq Ch^{\ell-k} | z |_{\ell,p} 
\qquad && \forall\ z \in W^{\ell,p}(I). \label{eq:estpih} 
\end{alignat}
\end{subequations}
In addition, for $z \in H^1(I)$ it holds that
\begin{equation} \label{eq:mvt}
\int_I (z - \pi^h z)_\rho \eta_{\rho} \drho = \sum_{j=1}^J
\int_{I_j} (z- \pi^h z)_\rho \eta_{\rho} \drho 
= 0 \qquad \forall\ \eta \in V^h.
\end{equation}
A natural semidiscrete finite element approximation for the system 
\eqref{eq:sys} is now defined as follows: 
Let $x^0_h \in \Vh$ be given. Then find
$x_h, y_h: I \times [0,T] \rightarrow \bR^d$ such that $x_h(\cdot,t), y_h(\cdot,t) \in \Vh$, $x_h(\cdot,0)= x^0_h$ and for $t \in [0,T]$ it holds that
\begin{subequations} \label{eq:sd}
\begin{align} 
& \int_I x_{h,t} \cdot \chi |x_{h,\rho}|^2 \drho - \int_I y_{h,\rho}  \cdot \chi_\rho \drho
=  \int_I F(x_{h,\rho},y_h,y_{h,\rho})y_h \cdot \chi \drho 
\quad \forall\ \chi \in \Vh, \label{eq:sda} \\
& \int_I y_h \cdot \eta |x_{h,\rho}|^2 \drho + \int_I x_{h,\rho} \cdot \eta_\rho \drho
= 0
\qquad \forall\ \eta \in \Vh. \label{eq:sdb}
\end{align}
\end{subequations}
It turns out that for our main result, Theorem~\ref{thm:main} below,
the initial data $x^0_h \in \Vh$ has to be chosen very carefully, see also
Remark~\ref{rem:x0}.
We observe that in the case of curve diffusion, $F=F_{cd}$, a semidiscrete 
analogue of \eqref{eq:dtE0} holds. Indeed, choosing 
$\chi = y_h$ in \eqref{eq:sda} and $\eta = x_{h,t}$ in
\eqref{eq:sdb} yields, on recalling \eqref{eq:antisymm}, that
\begin{align}\label{eq:stab}
\tfrac12 \ddt \int_I |x_{h,\rho}|^2 \drho & = \int_I x_{h,\rho} \cdot x_{h,t\rho} \drho
= - \int_I y_h \cdot x_{h,t}  |x_{h,\rho} |^2 \drho 
\nonumber \\ 
& = - \int_I |y_{h,\rho} |^2 \drho - \int_I F_1(x_{h,\rho},y_h,y_{h,\rho}) y_h \cdot y_h \drho \nonumber \\
& = - \int_I |y_{h,\rho} |^2 \drho - 2 \int_I y_{h,\rho}\cdot x_{h,\rho} |y_h|^2 \drho
- \int_I |y_h|^4|x_{h,\rho} |^2 \drho \nonumber \\ &
= - \int_I |y_{h,\rho} + |y_h|^2 x_{h,\rho}|^2 \drho \leq 0.
\end{align}

Our main result are the following optimal error bounds which are valid both for the case of curve diffusion and for the elastic flow.

\begin{theorem} \label{thm:main}
Let $F=F_{cd}$ or $F=F_{el}$, recall \eqref{eq:defF} and \eqref{eq:defF3}. 
Suppose that $x^0_h=\hat x^0_h$ where $\hat x^0_h \in \Vh$ solves
\begin{equation} \label{eq:defQh}
 \int_I \hat x^0_{h,\rho} \cdot \eta_\rho \drho 
+ \int_I \hat x^0_h  \cdot \eta \drho =\int_I \pi^h x_0 \cdot \eta \drho - \int_I \pi^h\Bigl[\frac{x_{0,\rho\rho}}{|x_{0,\rho}|^2}\Bigr] \cdot \eta \, | (\pi^h x_0)_\rho |^2 \drho 
 \qquad \forall\ \eta \in \Vh.
\end{equation}
Then
there exists $h_0>0$ such that for $0 < h \leq h_0$ the semidiscrete problem \eqref{eq:sd} has a unique solution $(x_h,y_h): I \times [0,T] \rightarrow \bR^d \times \bR^d$ and the
following error bounds hold:
\begin{subequations} \label{eq:main}
\begin{align}
\max_{t \in [0,T]} \| x(\cdot,t) - x_h(\cdot,t) \|_0^2 + h^2 \max_{t \in [0,T]} | x(\cdot,t) - x_h(\cdot,t) |_1^2 + \int_0^T \| x_t - x_{h,t} \|_0^2 \dt & \leq C h^4, \quad \label{eq:ebx} \\
\max_{t \in [0,T]} \| y(\cdot,t) - y_h(\cdot,t) \|_0^2 + h^2 \max_{t \in [0,T]}  | y(\cdot,t) - y_h(\cdot,t) |_1^2 & \leq Ch^4. \label{eq:erry}
\end{align}
\end{subequations}
\end{theorem}

\begin{remark} \label{rem:x0}
Note that the initial datum $y_h^0=y_h(\cdot,0)$ is not prescribed, 
but is determined by $x_h(\cdot,0)=x_h^0$ through the relation 
\eqref{eq:sdb}, i.e.\
\begin{equation*} 
\int_I y_h^0 \cdot \eta |x^0_{h,\rho}|^2 \drho 
+ \int_I x^0_{h,\rho} \cdot \eta_\rho \drho = 0
\qquad \forall\ \eta \in \Vh.
\end{equation*}
We will show in the Appendix, see Lemma~\ref{lem:Qh}, that 
with the choice $x^0_h = \hat x^0_h$ given by \eqref{eq:defQh} one has for sufficiently small $0<h \leq h_0$ that
$y^0_h$ exists uniquely and satisfies
\begin{equation} \label{eq:estinit}
 \| y(\cdot,0) - y_h(\cdot,0) \|_0 \leq C h^2,
\end{equation}  
which is crucial for obtaining optimal error estimates.
It does not seem straightforward to prove the bound \eqref{eq:estinit} for 
the more natural choice $x_h^0 = \pi^h x_0$.
\end{remark}

\setcounter{equation}{0}
\section{Proof of Theorem \ref{thm:main}} \label{sec:proof}

For $w \in H^1(I;\bR^d)$ we denote in what follows by $R_h w \in \Vh$ the solution to 
\begin{align} 
& \int_I (R_h w)_\rho \cdot \chi_\rho \drho 
+ \int_I F(x_\rho,y,(R_h w)_\rho) y \cdot \chi \drho 
+ \gamma \int_I R_h w \cdot \chi \drho \nonumber \\ & \qquad
= \int_I w_\rho \cdot \chi_\rho \drho + \int_I F(x_\rho,y,w_\rho) y \cdot \chi \drho + \gamma \int_I  w \cdot \chi \drho \quad \forall \chi \in \Vh. 
\label{eq:defrh}
\end{align}
Note that the operator $R_h$ depends on $x$ and $y$.
It is shown in  Lemma~\ref{lem:Rh} that $R_h w$ exists uniquely provided that $\gamma \geq 18 C_0^4+\frac{1}{2}$.
Let us decompose the errors $x-x_h$ and $y-y_h$ as
\begin{align*}
x-x_h& = (x-\pi^h x) + (\pi^h x-x_h)=: \xpihx+ \pihxxh,  \\
y-y_h & = (y-R_h y) + (R_h y - y_h)= :\yRhy + \Rhyyh.
\end{align*}
Lemma~\ref{lem:Rh} yields for $t \in [0,T]$ that
\begin{equation} \label{eq:estrh}
\| \yRhy \|_0 + h  | \yRhy |_1 + \| \yRhy_t \|_0 \leq C h^2,
\end{equation}  
while $\pihxxh, \Rhyyh \in \Vh$. As a result of the particular form of $R_h y$ we shall be able to prove the superconvergence property that $\pihxxh_\rho$ and $\Rhyyh_\rho$ are $\mathcal{O}(h^2)$, which in turn will be
crucial in order to derive optimal $L^2$--error bounds. Roughly speaking, our use of $\pi^h x$ and $R_h y$ can be viewed as a generalization of the use of the Ritz projection introduced by Wheeler in \cite{Wheeler73} for the heat equation.

Taking the scalar product of  \eqref{eq:sys1} with $\chi \in \Vh$, integrating over $I$ and using  \eqref{eq:defrh} we obtain
\begin{equation} \label{eq:weaka1} 
 \int_I x_t \cdot \chi |x_\rho|^2 \drho - \int_I (R_h y)_\rho \cdot \chi_\rho \drho = \int_I F(x_\rho,y,(R_h y)_\rho) y \cdot \chi \drho 
- \gamma \int_I \yRhy \cdot \chi \drho \quad \forall\ \chi \in \Vh. 
\end{equation}

Let us define
\begin{align} \label{eq:hatTh}
\hat T_h = \sup \Bigl\{ t \in [0,T] & :
(x_h,y_h) \text{ solves } (\ref{eq:sd}) \text{ on } [0,t],
\text{ with }
\int_0^t  |  x_{h,t} |_{1,\infty}^2 \du \leq 2 C_0\text{ and } 
\nonumber \\ & \quad
 \tfrac12 c_0 \leq  | x_{h,\rho} |
\ \text{and}\
| x_{h,\rho} | + | y_h |  \leq 2 C_0,
\ 0 \leq u \leq t \Bigr\}.
\end{align}

We will split the error analysis into a series of auxiliary results. To begin, let us derive two error relations. Taking the difference of \eqref{eq:weaka1} and \eqref{eq:sda} we obtain 
\begin{align}
& \int_I \pihxxh_t \cdot \chi  | x_{h,\rho}|^2 \drho - \int_I \Rhyyh_\rho \cdot \chi_\rho \drho 
= -\int_I \sigma_t \cdot \chi | x_{h,\rho} |^2 \drho +  \int_I x_t \cdot \chi \bigl[ | x_{h,\rho} |^2 - | x_\rho |^2 \bigr] \drho 
\nonumber \\ & \qquad 
+ \int_I \bigl( F(x_\rho,y,(R_h y)_\rho)y - F(x_{h,\rho},y_h, y_{h,\rho}) y_h \bigr) \cdot \chi \drho - \gamma \int_I \yRhy \cdot \chi \drho \qquad \forall\ \chi \in \Vh.
\label{eq:erra}
\end{align}

Furthermore, we infer from \eqref{eq:sys2}, \eqref{eq:sdb} and \eqref{eq:mvt}
that 
\begin{equation} \label{eq:errb}
 \int_I \Rhyyh \cdot \eta | x_{h,\rho} |^2 \drho + \int_I \pihxxh_\rho \cdot \eta_\rho \drho = - \int_I \yRhy \cdot \eta | x_{h,\rho} |^2 \drho + 
 \int_I y \cdot \eta \bigl[ | x_{h,\rho} |^2 - | x_\rho |^2 \bigr] \drho \quad \forall \eta \in \Vh.
\end{equation}

\begin{lemma}  \label{lem:l1}
It holds that
\[
\tfrac12 \ddt | \pihxxh |_1^2 \leq  
2\ddt \int_I ( y \cdot \pihxxh_\rho) x_\rho \cdot \xpihx \drho 
+ \epsilon \| \pihxxh_t \|_0^2 + C_\epsilon h^4 + C_\epsilon \| \Rhyyh \|_0^2 + C_\epsilon | \pihxxh |_1^2.
\]
\end{lemma}
\begin{proof}
Inserting $\eta=\pihxxh_t$ into \eqref{eq:errb} we infer with the help of integration by parts
\begin{align*}
 \tfrac12 \ddt | \pihxxh |_1^2 &
= - \int_I \yRhy \cdot \pihxxh_t \, | x_{h,\rho} |^2 \drho - \int_I \Rhyyh \cdot  \pihxxh_t \, | x_{h,\rho} |^2 \drho+ \int_I y \cdot \pihxxh_t \bigl[ | x_{h,\rho} |^2 - | (\pi^h x)_\rho |^2 \bigr] \drho \\ & \quad
+ \int_I  y \cdot \pihxxh_t | \xpihx_\rho |^2 \drho -2 \int_I  (y \cdot \pihxxh_t) x_\rho \cdot \xpihx_\rho \drho \\ & 
\leq C \bigl( \| \yRhy \|_0 + \| \Rhyyh \|_0 + | \pihxxh |_1 + | \xpihx |_1 | \xpihx |_{1,\infty} \bigr) \| \pihxxh_t \|_0 \\ & \quad 
+2 \int_I (y \cdot \pihxxh_{t\rho})  x_\rho \cdot \xpihx \drho 
+2 \int_I (y \cdot \pihxxh_t) x_{\rho\rho} \cdot \xpihx \drho 
+2 \int_I (y_\rho \cdot \pihxxh_t)  x_\rho \cdot \xpihx \drho \\ &
\leq C \bigl( \| \Rhyyh \|_0 + | \pihxxh |_1 + h^2 \bigr) \| \pihxxh_t \|_0 +2 \int_I (y \cdot \pihxxh_{t\rho})  x_\rho \cdot \xpihx \drho,
\end{align*}
where we also used  \eqref{eq:estrh}, \eqref{eq:hatTh} and \eqref{eq:estpih}.
Finally, estimating
\begin{align*}
& \int_I (y \cdot \pihxxh_{t\rho})  x_\rho \cdot \xpihx \drho \\ & \quad
=  \ddt \int_I ( y \cdot \pihxxh_\rho) x_\rho \cdot \xpihx  \drho - \int_I (y_t \cdot \pihxxh_\rho) x_\rho \cdot \xpihx \drho - \int_I ( y \cdot \pihxxh_\rho) x_{t,\rho} \cdot \xpihx  \drho
- \int_I ( y \cdot \pihxxh_\rho) x_\rho \cdot \xpihx_t  \drho \\ & \quad
\leq \ddt \int_I ( y \cdot \pihxxh_\rho) x_\rho \cdot \xpihx  \drho + C h^2 | \pihxxh |_1,
\end{align*}
by \eqref{eq:estpih}, we deduce the desired result with the help of Young's inequality.
\end{proof}

\begin{lemma} \label{lem:l3}
It holds that
\begin{displaymath}
 | \Rhyyh |_1^2 \leq \epsilon \| \pihxxh_t \|_0^2 + C_\epsilon \| \Rhyyh \|_0^2 + C h^4 + C | \pihxxh |_1^2. 
\end{displaymath}
\end{lemma}
\begin{proof}
Inserting $\chi= \Rhyyh$ into \eqref{eq:erra} and rearranging yields
\begin{align*}
| \Rhyyh |_1^2 & = \int_I \pihxxh_t \cdot \Rhyyh | x_{h,\rho}|^2 \drho + \int_I \sigma_t \cdot \Rhyyh | x_{h,\rho} |^2 \drho -  \int_I x_t \cdot \Rhyyh \bigl[ | x_{h,\rho} |^2 - | x_\rho |^2 \bigr] \drho 
\nonumber \\ & \quad
 - \int_I \bigl( F(x_\rho,y,(R_h y)_\rho)y - F(x_{h,\rho},y_h, y_{h,\rho}) y_h \bigr) \cdot \Rhyyh \drho + \gamma \int_I \yRhy \cdot \Rhyyh \drho 
=: \sum_{i=1}^5 \tilde T_i . 
\end{align*}
Clearly, \eqref{eq:hatTh}, \eqref{eq:estpih} and \eqref{eq:estrh} imply that 
\begin{equation} \label{eq:tildeT1}
| \tilde T_1 | + | \tilde T_2 | + | \tilde T_5 | \leq C \bigl( \Vert \pihxxh_t \Vert_0 + \Vert \sigma_t \Vert_0 + \Vert \yRhy \Vert_0 \bigr) \Vert \Rhyyh \Vert_0 
 \leq  \epsilon \| \pihxxh_t \|_0^2 + C_\epsilon \| \Rhyyh \|_0^2 + C h^4. 
\end{equation}
Next, let us write with the help of integration by parts
\begin{align} \label{eq:tildeT3}
\tilde T_3 & = 
  -\int_I x_t \cdot \Rhyyh \bigl[ | x_{h,\rho} |^2 - | (\pi^h x)_\rho |^2 \bigr] \drho 
- \int_I x_t \cdot \Rhyyh \, | \xpihx_\rho |^2 \drho
+ 2 \int_I (x_t \cdot \Rhyyh) \xpihx_\rho \cdot x_\rho \drho 
\nonumber \\ & 
= -\int_I  x_t \cdot \Rhyyh \bigl[ | x_{h,\rho} |^2 - | (\pi^h x)_\rho |^2 \bigr]  \drho - \int_I x_t \cdot \Rhyyh \, | \xpihx_\rho |^2 \drho 
- 2 \int_I (x_{t\rho} \cdot \Rhyyh) \xpihx  \cdot x_\rho \drho  
\nonumber\\ & \quad 
- 2 \int_I (x_t \cdot \Rhyyh_\rho) \xpihx  \cdot x_\rho \drho 
- 2 \int_I (x_t \cdot \Rhyyh) \xpihx \cdot x_{\rho \rho} \drho.
\end{align}
Hence, \eqref{eq:hatTh} and \eqref{eq:estpih} imply 
\begin{equation} \label{eq:tt1}
| \tilde T_3 |  \leq C \bigl( | \pihxxh |_1 + h^2 \bigr)  \| \Rhyyh \|_0 + C h^2 | \Rhyyh |_1 
 \leq \tfrac16 | \Rhyyh |_1^2 + C h^4 + C | \pihxxh |_1^2 + C \| \Rhyyh \|_0^2.
\end{equation}
In order to treat the remaining term $\tilde T_4$, we write 
\begin{align}
& F(x_\rho,y,(R_h y)_\rho)y - F(x_{h,\rho},y_h, y_{h,\rho}) y_h 
\nonumber \\ & \quad
= \left[ (F(x_\rho,y,(R_h y)_\rho)y - F((\pi^hx)_\rho,y, (R_h y)_\rho) y )  
 - (F(x_\rho,y,y_\rho)y - F((\pi^hx)_\rho,y, y_\rho) y ) \right] 
\nonumber \\ & \qquad
  +  \left[ F(x_\rho,y,y_\rho)y - F((\pi^h x)_\rho,y, y_\rho) y \right] 
  + \left[ F((\pi^h x)_\rho,y,(R_h y)_\rho)y -  F(x_{h,\rho},y_h, y_{h,\rho}) y_h \right] \nonumber \\ & \quad
  =: A_1+A_2 + A_3.  \label{eq:difF}
\end{align}
Using the intermediate value theorem we obtain that
\begin{equation} \label{eq:a1}
| A_1 | \leq C | \xpihx_\rho | \bigl( | \xpihx_\rho |+  | \yRhy_\rho | \bigr),
\end{equation}
so that \eqref{eq:estpih} and \eqref{eq:estrh} yield
\begin{displaymath}
| \int_I A_1 \cdot \Rhyyh \drho |  \leq C | \xpihx |_{1,\infty}  \bigl( | \yRhy |_1+ | \xpihx |_1 \bigr)  \Vert \Rhyyh \Vert_0 \leq C h^2 \Vert \Rhyyh \Vert_0 \leq Ch^4 + C \Vert \Rhyyh \Vert_0^2.
\end{displaymath}
Next, let us write
\begin{equation} \label{eq:a2}
A_2 = - G \xpihx_\rho + r, \quad \mbox{ where } G := \nabla_a [F(a,y,y_\rho)y]_{a=x_\rho} \mbox{ and } 
| r | \leq | \xpihx_\rho |^2,
\end{equation}
from which we deduce with the help of integration by parts and \eqref{eq:estpih} that 
\begin{align*}
| \int_I A_2 \cdot \Rhyyh \drho | & \leq | - \int_I G \xpihx_\rho \cdot \Rhyyh \drho | + | \int_I r \cdot \Rhyyh \drho | \leq | \int_I (G_\rho \xpihx \cdot \Rhyyh + G \xpihx \cdot \Rhyyh_\rho \drho |
+ C h^2 \Vert \Rhyyh \Vert_0 \\
& \leq C h^2 \Vert \Rhyyh \Vert_1 \leq \tfrac16 | \Rhyyh |_1^2 + C h^4 + \Vert \Rhyyh \Vert_0^2.
\end{align*}
Finally, since $| x_{h,\rho}| + |y_h |  \leq 2 C_0$ and $F$ depends linearly on $c$ we have
\begin{equation} \label{eq:a3}
| A_3 | \leq C \bigl( | \pihxxh_\rho | + | \yRhy | + | \Rhyyh | + | \Rhyyh_\rho | \bigr),
\end{equation}
so that \eqref{eq:estrh} implies
\begin{displaymath}
| \int_I A_3 \cdot \Rhyyh \drho |  \leq C \bigl( | \pihxxh |_1 + \Vert \yRhy \Vert_0 + \Vert \Rhyyh \Vert_0 + | \Rhyyh |_1 \bigr) \Vert \Rhyyh \Vert_0 \leq 
\tfrac16 | \Rhyyh |_1^2
+ C | \pihxxh |_1^2 + C \Vert \Rhyyh \Vert_0^2+C h^4.
\end{displaymath}
In conclusion 
\begin{equation} \label{eq:tildeT4}
| \tilde T_4 | \leq \tfrac13 | \Rhyyh |_1^2
+ C | \pihxxh |_1^2 + C \Vert \Rhyyh \Vert_0^2+C h^4.
\end{equation}
Combining \eqref{eq:tildeT1}, \eqref{eq:tt1} and \eqref{eq:tildeT4} yields 
the desired result.
\end{proof}

\begin{lemma}  \label{lem:l2}
It holds that
\[
\tfrac12 \ddt \int_I | \Rhyyh |^2 | x_{h,\rho} |^2 \drho + \int_I \pihxxh_{t\rho} \cdot \Rhyyh_\rho \drho 
\leq \epsilon \| \pihxxh_t \|_0^2 + C_\epsilon(1+ | x_{h,t} |_{1,\infty}^2) \bigl(  h^4 +  \| \Rhyyh \|_0^2 +  | \pihxxh |_1^2 \bigr).
\]
\end{lemma}
\begin{proof} 
Direct calculation shows that
\begin{equation} \label{eq:dtterm}
\tfrac12 \ddt \int_I | \Rhyyh |^2 | \, x_{h,\rho} |^2 \drho
= \int_I \Rhyyh \cdot \Rhyyh_t \, | x_{h,\rho} |^2 \drho
+ \int_I | \Rhyyh |^2 \, x_{h,\rho} \cdot x_{h,t\rho} \drho .
\end{equation}
Moreover, differentiating \eqref{eq:errb} with respect to $t$ we obtain
\begin{align}
& \int_I \Rhyyh_t  \cdot \eta | x_{h,\rho} |^2 \drho 
+ 2 \int_I (\Rhyyh \cdot \eta) x_{h,\rho}  \cdot x_{h,\rho t} \drho
+ \int_I \pihxxh_{t\rho}  \cdot \eta_\rho \drho  \nonumber \\ & \quad
 = - \int_I \yRhy_t \cdot \eta | x_{h,\rho}|^2 
 - 2 \int_I ( \yRhy \cdot \eta) x_{h,\rho} \cdot x_{h,\rho t} \drho 
 \nonumber \\ & \qquad 
+ \int_I y_t \cdot \eta  \bigl[ | x_{h,\rho} |^2 - | x_\rho |^2 \bigr] \drho
    + 2 \int_I y \cdot \eta  \bigl[ x_{h,\rho} \cdot x_{h,\rho t} - x_\rho \cdot x_{\rho t} \bigr] \drho.  \label{eq:errbt}
\end{align}
Combining \eqref{eq:dtterm} and \eqref{eq:errbt} with $\eta=\Rhyyh$ 
yields that
\begin{align*}
& \tfrac12 \ddt \int_I | \Rhyyh |^2 | x_{h,\rho} |^2 \drho + \int_I \pihxxh_{t\rho}  \cdot \Rhyyh_\rho \drho 
\\ & \quad
=- \int_I  \yRhy_t \cdot \Rhyyh | x_{h,\rho} |^2 -  \int_I | \Rhyyh |^2 x_{h,\rho} \cdot x_{h,\rho t} \drho 
- 2 \int_I ( \yRhy \cdot \Rhyyh ) x_{h,\rho} \cdot x_{h,\rho t} \drho 
\\ & \qquad 
+ \int_I (y_t \cdot \Rhyyh) ( | x_{h,\rho} |^2 - | x_\rho |^2 ) \drho
+ 2 \int_I y \cdot \Rhyyh \bigl[ x_{h,\rho} \cdot x_{h,\rho t} - x_\rho \cdot x_{\rho t} \bigr] \drho =: \sum_{i=1}^5 S_i.
\end{align*}
To begin, we deduce from \eqref{eq:estrh} that
\begin{align}
| S_1 | + | S_2 | + | S_3 | &\leq C \| \yRhy_t\|_0 \| \Rhyyh \|_0 + C \bigl(  \| \Rhyyh \|_0^2 + h^2 \| \Rhyyh \|_0 \bigr)  | x_{h,t} |_{1,\infty}  \nonumber \\
& \leq C h^4 +  C \bigl( 1+ |x_{h,t}|_{1,\infty}^2 \bigr) \| \Rhyyh \|_0^2.  \label{eq:s23}
\end{align}
Next, similarly to \eqref{eq:tildeT3} and \eqref{eq:tt1}, 
we infer 
\begin{align} \label{eq:s4} 
|  S_4 | & \leq C \| \Rhyyh \|_0 | \pihxxh |_1 + Ch^2 \| \Rhyyh  \|_1 \leq C h^4 + C | \pihxxh|_1^2 + C \| \Rhyyh \|_1^2 \nonumber \\
& \leq  \tfrac12 {\epsilon} \| \pihxxh_t \|_0^2 + C_\epsilon  \| \Rhyyh \|_0^2 + C h^4 + C  | \pihxxh |_1^2,
\end{align} 
where we also used Lemma~\ref{lem:l3}. Finally,
\begin{align*}
S_5 & = - 2 \int_I (y \cdot \Rhyyh) x_\rho \cdot (\xpihx_{t\rho} + \pihxxh_{t\rho}) \drho 
 - 2 \int_I (y \cdot \Rhyyh) \xpihx_\rho \cdot x_{h,\rho t} \drho 
\\ & \quad 
- 2 \int_I (y \cdot \Rhyyh) \pihxxh_\rho \cdot x_{h,\rho t} \drho 
=: S_{5,1} +  S_{5,2} +  S_{5,3}.
\end{align*}
Using integration by parts we deduce with \eqref{eq:estpih} that
\begin{displaymath}
| S_{5,1} | \leq C \| \Rhyyh \|_1 \bigl( \| \xpihx_t \|_0+ \| \pihxxh_t \|_0 \bigr)  \leq \tfrac14 {\epsilon} \| \pihxxh_t \|_0^2 +C_\epsilon \| \Rhyyh \|_1^2 
+ Ch^4.
\end{displaymath}
On recalling \eqref{eq:mvt} we may write with $c_j = h_j^{-1} \int_{I_j} y \cdot \Rhyyh \drho$
\begin{displaymath}
 S_{5,2}  = -2 \sum_{j=1}^J \int_{I_j} \left[ y \cdot \Rhyyh - c_j  \right] \xpihx_\rho \cdot x_{h,\rho t} \drho 
\leq C h |y \cdot \Rhyyh|_1 |\xpihx|_1 |x_{h, t}|_{1,\infty}
\leq  C h^2  | x_{h,t} |_{1,\infty}\| \Rhyyh \|_1.
\end{displaymath}
Since $| S_{5,3} | \leq C | x_{h,t} |_{1,\infty} \| \Rhyyh \|_0 | \pihxxh |_1$ we thus have
\begin{align} \label{eq:s5}
| S_5 | & \leq \tfrac14 {\epsilon} \| \pihxxh_t \|_0^2 + C_\epsilon \| \Rhyyh \|_1^2 
+ C (1+ | x_{h,t} |_{1,\infty}^2) h^4 
+ C_\epsilon  | x_{h,t} |_{1,\infty}^2  | \pihxxh |_1^2 \nonumber \\
& \leq  \tfrac12 {\epsilon} \| \pihxxh_t \|_0^2 + C_\epsilon \| \Rhyyh \|_0^2 + C_\epsilon (1+ | x_{h,t} |_{1,\infty}^2) h^4
+ C_\epsilon (1+ | x_{h,t} |_{1,\infty}^2) | \pihxxh |_1^2,
\end{align}
where we again used Lemma~\ref{lem:l3}.
Combining the estimates in \eqref{eq:s23}, 
\eqref{eq:s4} and \eqref{eq:s5} yields the desired result.
\end{proof}

\begin{lemma}  \label{lem:l4}
Let $G$ be as in \eqref{eq:a2}. Then we have 
\begin{align*}
\tfrac18 c_0^2 \| \pihxxh_t \|_0^2 - \int_I \Rhyyh_\rho \cdot \pihxxh_{t\rho} \drho  & \leq  \ddt \int_I \left[ 2(x_t \cdot \pihxxh_\rho) \xpihx \cdot x_\rho + 
G \xpihx \cdot \pihxxh_{\rho} \right] \drho \\
& \qquad + Ch^4 + C \bigl( 1+ \Vert x_{tt} \Vert_0^2 \bigr)  | \pihxxh |_1^2  + C \| \Rhyyh \|_0^2.
\end{align*}
\end{lemma}
\begin{proof}
Choosing $\chi= \pihxxh_t$ in \eqref{eq:erra} yields together with \eqref{eq:hatTh} 
\begin{align}
& \tfrac14 c_0^2 \| \pihxxh_t \|_0^2 - \int_I \Rhyyh_\rho \cdot \pihxxh_{t\rho} \drho 
 \leq  - \int_I \xpihx_t \cdot \pihxxh_t | x_{h,\rho} |^2 \drho 
+ \int_I x_t \cdot \pihxxh_t \bigl[ | x_{h,\rho} |^2 - | x_\rho |^2 \bigr] \drho \nonumber \\  & \quad
+  \int_I \bigl( F(x_\rho,y,(R_h y)_\rho)y - F(x_{h,\rho},y_h, y_{h,\rho}) y_h \bigr) \cdot \pihxxh_t \drho - \gamma \int_I \yRhy \cdot \pihxxh_t \drho 
 =: \sum_{i=1}^{4} T_i. \label{eq:err1}
\end{align}
Clearly, it follows from \eqref{eq:estpih}, \eqref{eq:hatTh} and \eqref{eq:estrh} that
\begin{equation} \label{eq:t1}
| T_1 | + | T_4 |  \leq C \bigl(  \| \xpihx_t \|_0 + \Vert \yRhy \Vert_0 \bigr)  \| \pihxxh_t \|_0 \leq \tfrac{1}{24} \| \pihxxh_t \|_0^2 + C h^4.
\end{equation}
Next, arguing in a similar way as for $\tilde T_3$ in \eqref{eq:tildeT3} 
we write
\begin{align*}
T_2 & = \int_I  x_t \cdot \pihxxh_t \bigl[ | x_{h,\rho} |^2 - | (\pi^h x)_\rho |^2 \bigr]  \drho + \int_I x_t \cdot \pihxxh_t \, | \xpihx_\rho |^2 \drho 
+ 2 \int_I (x_{t\rho} \cdot \pihxxh_t) \xpihx  \cdot x_\rho \drho  
\nonumber\\ & \quad 
+ 2 \int_I (x_t \cdot \pihxxh_{t\rho}) \xpihx  \cdot x_\rho \drho 
+ 2 \int_I (x_t \cdot \pihxxh_t) \xpihx \cdot x_{\rho \rho} \drho
=: \sum_{j=1}^5 T_{2,j}.\nonumber
\end{align*} 
Clearly,
\begin{displaymath}
\sum_{j\in\{1,2,3,5\}} | T_{2,j} | \leq  C \| \pihxxh_t \|_0 | \pihxxh |_1 + C h^2 \| \pihxxh_t \|_0, 
\end{displaymath}
while we write for the remaining term
\begin{align*}
T_{2,4} & = 2 \ddt \int_I (x_t \cdot \pihxxh_\rho) \xpihx \cdot x_\rho \drho -2 \int_I (x_{tt} \cdot \pihxxh_{\rho}) \sigma \cdot x_\rho  +  (x_t \cdot \pihxxh_{\rho}) \sigma_t \cdot x_\rho  +  (x_t \cdot \pihxxh_{\rho}) \sigma \cdot x_{t\rho} \drho \\
& \leq 2 \ddt \int_I (x_t \cdot \pihxxh_\rho) \xpihx \cdot x_\rho \drho  + C h^2 \bigl( 1+  \Vert x_{tt} \Vert_0 \bigr)   | \pihxxh |_1,
\end{align*}
so that
\begin{align} 
T_2  &\leq 2 \ddt \int_I (x_t \cdot \pihxxh_\rho) \xpihx \cdot x_\rho \drho 
+ C \| \pihxxh_t \|_0 | \pihxxh |_1 +  C h^2 \| \pihxxh_t \|_0 + C h^2 \bigl( 1 + \Vert x_{tt} \Vert_0 \bigr)  | \pihxxh |_1  
\nonumber \\ 
& \leq 2 \ddt \int_I (x_t \cdot \pihxxh_\rho) \xpihx \cdot x_\rho \drho  + \tfrac{1}{24} \| \pihxxh_t \|_0^2 + C  h^4 + C \bigl( 1 + \Vert x_{tt} \Vert_0^2 \bigr)  | \pihxxh |_1^2. \label{eq:t2}
\end{align} 
Finally, we write $  T_3  = \sum_{j=1}^3 \int_I A_j \cdot \pihxxh_t \drho$
with $A_1,A_2, A_3$ as in \eqref{eq:difF}. We obtain with the help of \eqref{eq:a1} and \eqref{eq:a3} that 
\begin{displaymath}
\left| \int_I (A_1+A_3) \cdot \pihxxh_t \drho \right| \leq 
C \bigl(  | \xpihx |_{1,\infty} | \xpihx |_1 +
 | \xpihx |_{1,\infty} | \yRhy |_1 +| \pihxxh |_1 + \Vert \yRhy \Vert_0 + \Vert \Rhyyh \Vert_0 + | \Rhyyh |_1 \bigr)
 \Vert \pihxxh_t \Vert_0,
\end{displaymath}
while \eqref{eq:a2} yields
\begin{align*}
 \int_I A_2 \cdot \pihxxh_t \drho &=- \int_I G \xpihx_\rho \cdot \pihxxh_t \drho + \int_I r \cdot \pihxxh_t \drho \\
& = \int_I G_\rho \xpihx \cdot \pihxxh_t \drho + \int_I G \xpihx \cdot \pihxxh_{t,\rho}  \drho + \int_I r \cdot \pihxxh_t \drho \\
& = \int_I G_\rho \xpihx \cdot \pihxxh_t \drho + \ddt \int_I G \xpihx \cdot \pihxxh_{\rho}  \drho - \int_I G_t \xpihx \cdot \pihxxh_\rho \drho 
- \int_I G \xpihx_t \cdot \pihxxh_\rho \drho 
+ \int_I r \cdot \pihxxh_t \drho \\
& \leq \ddt \int_I G \xpihx \cdot \pihxxh_{\rho}  \drho + C h^2 \bigl(  \Vert \pihxxh_t \Vert_0 + |  \pihxxh |_1 \bigr).
\end{align*}
Combining the above estimates we obtain with the help of Lemma~\ref{lem:l3}
\begin{align}
T_3 &  \leq \ddt \int_I G \xpihx \cdot \pihxxh_{\rho}  \drho  + \tfrac{1}{32} \Vert \pihxxh_t \Vert_0^2 +C  | \pihxxh |_1^2 + C  \Vert \Rhyyh \Vert_0^2    +  C | \Rhyyh |_1^2 + C h^4  \nonumber \\
& \leq \ddt \int_I G \xpihx \cdot \pihxxh_{\rho}  \drho + \tfrac{1}{24} \Vert \pihxxh_t \Vert_0^2 + C  | \pihxxh |_1^2 + C  \Vert \Rhyyh \Vert_0^2 + C h^4.  \label{eq:t3}
\end{align}
If we insert \eqref{eq:t1}, \eqref{eq:t2}, \eqref{eq:t3} into  \eqref{eq:err1}
we obtain the desired result. 
\end{proof}

If we combine Lemmas~\ref{lem:l1}, \ref{lem:l2} and \ref{lem:l4}, 
we obtain after choosing $\epsilon$ sufficiently small
\begin{align*}
&
\tfrac1{16} c_0^2 \| \pihxxh_t \|_0^2 
+ \tfrac12 \ddt | \pihxxh |_1^2 
+ \tfrac12 \ddt \int_I | \Rhyyh |^2 | x_{h,\rho} |^2 \drho \\ & \qquad
 \leq   \varphi'(t)  + C \bigl( 1 + | x_{h,t} |_{1,\infty}^2 +  \Vert x_{tt} \Vert_0^2 \bigr) \bigl( h^4 +  | \pihxxh  |_1^2 +  \| \Rhyyh \|_0^2 \bigr),
\end{align*}
where 
\begin{displaymath}
\varphi(t) = 2  \int_I ( y \cdot \pihxxh_\rho) x_\rho \cdot \xpihx \drho + 2 \int_I (x_t \cdot \pihxxh_\rho) \xpihx \cdot x_\rho \drho + 
 \int_I G \xpihx \cdot \pihxxh_{\rho}  \drho
 \end{displaymath}
satisfies $|\varphi(t)| \leq Ch^2 | \pihxxh(t) |_1$. Observing that
\begin{displaymath}
| \pihxxh(0) |_1^2  + \Vert \Rhyyh(0) \Vert_0^2 \leq | \pi^h x_0 - \hat x_0^h |_1^2 + 2 \Vert R_h y(0) - y(0)\Vert_0^2 + 2 \Vert y(0) - y_h(0) \Vert_0^2 \leq C h^4
\end{displaymath}
in view of \eqref{eq:Appestqh}, \eqref{eq:Appestrha} and \eqref{eq:Appestyh}, 
we obtain after integration with respect to $u \in (0,t), 0 < t < \hat T_h$ 
\begin{align*}
& \tfrac1{16} c_0^2 \int_0^t \| \pihxxh_t \|_0^2 \du 
+ \tfrac12 | \pihxxh(t) |_1^2 
+ \tfrac1{8} c_0^2 \| \Rhyyh(t) \|_0^2 \\ 
& \quad \leq Ch^4 + Ch^2 | \pihxxh(t) |_1 + C \int_0^t (1+ | x_{h,t} |_{1,\infty}^2 +  \Vert x_{tt} \Vert_0^2  ) \bigl( h^4 +  | \pihxxh |_1^2 +  \| \Rhyyh \|_0^2 \bigr) \du.
\end{align*}
Absorbing $C h^2 | e^h(t) |_1$ into the left hand side, and recalling
\eqref{eq:hatTh} as well as our assumption $x_{tt} \in L^2(0,T; L^2(I,\mathbb R^d))$, we deduce that
\begin{equation*} 
\int_0^t \| \pihxxh_t \|_0^2 \du 
+ | \pihxxh(t) |_1^2 +  \| \Rhyyh(t) \|_0^2 
\leq C h^4 + C \int_0^t (1+ | x_{h,t} |_{1,\infty}^2 +  \Vert x_{tt} \Vert_0^2 ) \bigl( | \pihxxh |_1^2 
+ \| \Rhyyh \|_0^2 \bigr) \du . 
\end{equation*}
With the help of Gronwall's lemma and \eqref{eq:hatTh} we infer that
\begin{equation} \label{eq:errest}
\int_0^t \| \pihxxh_t\|_0^2 \du + \sup_{0 \leq t < \hat T_h} \bigl(  | \pihxxh(t) |_1^2 +  \| \Rhyyh(t) \|_0^2 \bigr) \leq C h^4, \quad 0 \leq t < \hat T_h.
\end{equation}
We are now in a position to prove that $\hat T_h =T$. Suppose that $\hat T_h <T$. 
In view of \eqref{eq:inverse} it holds that
\begin{displaymath}
| x_{h,t} |_{1,\infty} \leq | \pihxxh_t |_{1,\infty} + | \xpihx_t  |_{1,\infty} + | x_t |_{1,\infty} \leq ch^{-\frac{3}{2}} \Vert \pihxxh_t \Vert_0 + c h + | x_t |_{1,\infty}
\end{displaymath}
and so we deduce with the help of \eqref{eq:errest} that
\begin{displaymath}
\int_0^{\hat T_h} | x_{h,t} |_{1,\infty}^2 \dt \leq \tfrac{4}{3} \int_0^T | x_t |_{1,\infty}^2 \dt + C h^{-3} \Vert \pihxxh_t \Vert_0^2 \dt + C h^2 \leq \tfrac{4}{3} C_0 + Ch \leq \tfrac{3}{2} C_0,
\end{displaymath}
provided that $0<h\leq h_0$ and $h_0$ is sufficiently small. In a similar way one shows that $ | x_{h,\rho}(\hat T_h)| \geq \frac{2}{3} c_0$ as well as 
$| x_{h,\rho}(\hat T_h) | + | y_h(\hat T_h) |  \leq \frac{3}{2} C_0$. Continuing the discrete solution beyond $\hat T_h$ then yields a contradiction to
the definition of $\hat T_h$.
The bounds \eqref{eq:main} now follow from \eqref{eq:errest}, 
\eqref{eq:estpih} and \eqref{eq:estrh}, where we also use
\eqref{eq:inverse} to control $\sup_{0 \leq t \leq T} | \Rhyyh(t) |_1$.
This concludes the proof of Theorem~\ref{thm:main}.

\setcounter{equation}{0}
\section{Fully discrete approximation} \label{sec:fd}

\subsection{Curve diffusion}
A fully discrete approximation of \eqref{eq:sd} in the case $F=F_{cd}$
is given as follows.

For $m \geq 0$, given $(x^m_h, y^m_h) \in \Vh \times \Vh$ find 
$(x^{m+1}_h, y^{m+1}_h) \in \Vh \times \Vh$ such that
\begin{subequations} \label{eq:fd}
\begin{align}
& \int_I \frac{x^{m+1}_h-x^m_h}{\Delta t} \cdot \chi |x^m_{h,\rho}|^2 \drho 
- \int_I y^{m+1}_{h,\rho} \cdot \chi_\rho \drho 
= 2 \int_I (y^{m+1}_{h,\rho} \cdot x^m_{h,\rho}) y^m_h \cdot \chi \drho
\nonumber \\ & \qquad
+ \int_I | x^m_{h,\rho} |^2 (y^m_h \cdot y^{m+1}_h) y^m_h \cdot \chi \drho
+  \int_I F_2(x^m_{h,\rho},y^m_h,y^m_{h,\rho})y^{m+1}_h \cdot \chi \drho 
\qquad \forall\ \chi \in \Vh, \label{eq:fda} \\
& \int_I y^{m+1}_h \cdot \eta |x^m_{h,\rho}|^2 \drho 
+ \int_I x^{m+1}_{h,\rho} \cdot \eta_\rho \drho = 0
\qquad \forall\ \eta \in \Vh. \label{eq:fdb}
\end{align}
\end{subequations}
The choice of semi-implicit treatment for the fully discrete approximation
of the term on the right hand side of \eqref{eq:sda}, recall \eqref{eq:defF},
was guided by the aim to obtain a linear scheme that is unconditionally stable.
In fact, we have the following existence, uniqueness and stability result.

\begin{theorem} \label{thm:exsstab}
Assume that $|x^m_{h,\rho}| > 0$ in I.
Then there exists a unique solution 
$(x^{m+1}_h, y^{m+1}_h) \in \Vh \times \Vh$ to \eqref{eq:fd}.
Moreover, any solution to \eqref{eq:fd} satisfies the stability bound
\begin{equation} \label{eq:fdstab}
\tfrac12 \int_I |x^{m+1}_{h,\rho}|^2 \drho 
+ 
\Delta t \int_I | y^{m+1}_{h,\rho} + (y^m_h \cdot y^{m+1}_h) x^m_{h,\rho} |^2
\drho 
\leq \tfrac12 \int_I |x^m_{h,\rho}|^2 \drho.
\end{equation}
\end{theorem}
\begin{proof}
As \eqref{eq:fd} represents a square linear system, it is sufficient to prove
that only the zero solution solves the homogeneous system. 
Let $(X_h,Y_h) \in \Vh \times \Vh$ be such that
\begin{subequations} 
\begin{align}
\int_I \frac{X_h}{\Delta t} \cdot \chi |x^m_{h,\rho}|^2 \drho 
- \int_I Y_{h,\rho} \cdot \chi_\rho \drho   
& = 2 \int_I (Y_{h,\rho} \cdot x^m_{h,\rho}) y^m_h \cdot \chi \drho
+ \int_I | x^m_{h,\rho} |^2 (y^m_h \cdot Y_h) y^m_h \cdot \chi \drho
\nonumber \\ & \qquad
+ \int_I F_2(x^m_{h,\rho},y^m_h,y^m_{h,\rho}) Y_h \cdot \chi \drho 
\qquad \forall\ \chi \in \Vh, \label{eq:fd0a} \\
\int_I Y_h \cdot \eta |x^m_{h,\rho}|^2 \drho 
+ \int_I X_{h,\rho} \cdot \eta_\rho \drho & = 0
\qquad \forall\ \eta \in \Vh. \label{eq:fd0b}
\end{align}
\end{subequations}
Choosing $\chi = -\Delta t Y_h$ in \eqref{eq:fd0a}, $\eta = X_h$ in 
\eqref{eq:fd0b} and summing the two gives, on recalling 
\eqref{eq:antisymm}, that
\begin{align*}
0 & = \int_I |X_{h,\rho}|^2 \drho
+ \Delta t \int_I |Y_{h,\rho}|^2 \drho 
+ 2 \Delta t \int_I (Y_{h,\rho} \cdot x^m_{h,\rho}) y^m_h \cdot Y_h \drho
+ \Delta t \int_I | x^m_{h,\rho} |^2 (y^m_h \cdot Y_h)^2 \drho \\ & 
= \int_I |X_{h,\rho}|^2 \drho +
\Delta t \int_I | Y_{h,\rho} + (y^m_h \cdot Y_h) x^m_{h,\rho} |^2 \drho .
\end{align*}
It follows that $X_{h,\rho} = 0$ in $I$, and so first
\eqref{eq:fd0b} implies that $Y_h=0$, and then \eqref{eq:fd0a} implies that
$X_h=0$. Hence we have shown that there exists a unique solution to
\eqref{eq:fd}.

In order to prove \eqref{eq:fdstab}, we choose 
$\chi = y^{m+1}_h$ in \eqref{eq:fda} and $\eta = x^{m+1}_h - x^m_h$ in
\eqref{eq:fdb} to yield that 
\begin{align*}
& \tfrac12 \int_I |x^{m+1}_{h,\rho}|^2 \drho 
- \tfrac12 \int_I |x^m_{h,\rho}|^2 \drho 
\leq \int_I x^{m+1}_{h,\rho} \cdot (x^{m+1}_{h,\rho} - x^m_{h,\rho}) \drho
\nonumber \\ & \quad
= - \Delta t \int_I |y^{m+1}_{h,\rho}|^2 \drho 
- 2 \Delta t \int_I (y^{m+1}_{h,\rho} \cdot x^m_{h,\rho}) 
(y^m_h \cdot y^{m+1}_h) \drho
- \Delta t \int_I | x^m_{h,\rho} |^2 (y^m_h \cdot y^{m+1}_h)^2 \drho
\nonumber \\ & \quad
= 
- \Delta t \int_I | y^{m+1}_{h,\rho} + (y^m_h \cdot y^{m+1}_h) x^m_{h,\rho} |^2
\drho
\leq 0 , 
\end{align*}
where we have used once again \eqref{eq:antisymm}.
\end{proof}

Observe that \eqref{eq:fdstab} is a fully discrete analogue of \eqref{eq:stab},
recall also \eqref{eq:dtE0}. 
We note that a discrete analogue of the property 
$\ddt \int_I |x_\rho| \drho \leq 0$, is much harder to prove. 
For the simpler case of curve shortening flow, such a discrete analogue is 
shown in \cite[Lemma~4.1.3]{BanschDGP23} for the scheme originally proposed in
\cite{DeckelnickD95}. But at present it is not clear if these techniques can be
generalized from the second order flow to the fourth order problem studied
here. Nevertheless, we remark that in 
all our numerical experiments, both $\int_I |x^m_{h,\rho}|^2 \drho$
and $\int_I |x^m_{h,\rho}| \drho$ are monotonically decreasing.

\subsection{Elastic flow}
A fully discrete  approximation  of \eqref{eq:sd} in the case $F=F_{el}$
is given as follows.

For $m \geq 0$, given $(x^m_h, y^m_h) \in \Vh \times \Vh$ find 
$(x^{m+1}_h, y^{m+1}_h) \in \Vh \times \Vh$ such that
\begin{subequations} \label{eq:fdel}
\begin{align}
& \int_I \frac{x^{m+1}_h-x^m_h}{\Delta t} \cdot \chi |x^m_{h,\rho}|^2 \drho 
- \int_I y^{m+1}_{h,\rho} \cdot \chi_\rho \drho \nonumber \\ & \quad
= 2 \int_I (y^{m+1}_{h,\rho} \cdot x^m_{h,\rho}) y^m_h \cdot \chi \drho
+ \int_I | x^m_{h,\rho} |^2 (y^m_h \cdot y^{m+1}_h) y^m_h \cdot \chi \drho
\nonumber \\ & \qquad
+ \int_I F_2(x^m_{h,\rho},y^m_h,y^m_{h,\rho}) y^{m+1}_h \cdot \chi \drho + \int_I F_3(x^m_{h,\rho},y^m_h) y^m_h \cdot \chi \drho
\qquad \forall\ \chi \in \Vh, \label{eq:fdela} \\
& \int_I y^{m+1}_h \cdot \eta |x^m_{h,\rho}|^2 \drho 
+ \int_I x^{m+1}_{h,\rho} \cdot \eta_\rho \drho = 0
\qquad \forall\ \eta \in \Vh. \label{eq:fdelb}
\end{align}
\end{subequations}
Existence and uniqueness of a solution to \eqref{eq:fdel} can be shown as in 
Theorem~\ref{thm:exsstab}. However, as mentioned previously, solutions 
to \eqref{eq:xt2el} in general do not satisfy the Dirichlet energy
estimate \eqref{eq:dtE}, and so a stability bound of the form 
\eqref{eq:fdstab} does not hold for solutions of \eqref{eq:fdel}. 

\setcounter{equation}{0}
\section{Numerical results} \label{sec:nr}

We implemented \eqref{eq:fd} within the
finite element toolbox Alberta, \cite{Alberta}, using
the sparse factorization package UMFPACK, see \cite{Davis04},
for the solution of the linear systems of equations arising at each time level.

For all our numerical simulations we use a uniform partitioning of $[0,1]$,
so that $q_j = jh$, $j=0,\ldots,J$, with $h = \frac 1J$. Unless otherwise
stated, we use $J=512$, $\Delta t = 10^{-4}$ and choose $x_h^0 = \pi^h x_0$.
Given $x^0_h \in \Vh$, for the initial data $y^0_h \in \Vh$ we always choose
the solution of 
\begin{equation*}
\int_I y^{0}_h \cdot \eta |x^0_{h,\rho}|^2 \drho 
+ \int_I x^{0}_{h,\rho} \cdot \eta_\rho \drho = 0
\qquad \forall\ \eta \in \Vh,
\end{equation*}
compare with \eqref{eq:fdb}.
For our computed solutions we will often monitor the ratio
\begin{equation} \label{eq:ratio}
\ratio^m = \dfrac{\max_{j=1,\ldots,J} |x_h^m(q_j) - x_h^m(q_{j-1})|}
{\min_{j=1,\ldots, J} |x_h^m(q_j) - x_h^m(q_{j-1})|}
\end{equation}
between the lengths of the longest and shortest element. Clearly
$\ratio^m\geq1$, with equality if and only if the curve is equidistributed.
Moreover, at times we will also be interested in a possible blow-up in 
curvature for the fourth order evolution we approximate. 
To this end, given $x^{m}_h \in \Vh$, we introduce the discrete 
curvature vector $\varkappa^{m}_h \in \Vh$ such that
\begin{equation*} 
\int_I \pi^h [ \varkappa^{m}_{h} \cdot \eta_h ]\, |x^{m}_{h,\rho}| \drho
+ \int_I \frac{x^{m}_{h,\rho} \cdot \eta_{h,\rho}}{|x^{m}_{h,\rho}|} \drho =0
\qquad \forall\ \eta_h \in \Vh.
\end{equation*}
In practice we will then monitor the quantity
\begin{equation} \label{eq:kappainv}
K^m_\infty = \max_{j=1,\ldots,J} |\varkappa^m_h(q_j)|
\end{equation}
as an approximation to the maximal value of 
$|\varkappa| = \frac{|\tau_\rho|}{|x_\rho|}$. 

\newcommand{\errorxL}{\| x -  x_h\|_0}
\newcommand{\errorxH}{\| x -  x_h\|_1}
\newcommand{\erroryL}{\| y -  y_h\|_0}
\newcommand{\erroryH}{\| y -  y_h\|_1}
\subsection{Curve diffusion}

\begin{table}[!b]
\center
\begin{tabular}{|r|c|c|c|c|c|c|c|c|}
\hline
$J$ & $\errorxL$ & EOC  & $\errorxH$ & EOC  & $\erroryL$ & EOC  & $\erroryH$ & EOC \\ \hline
32  &3.5628e-03& ---&7.2517e-01& ---&3.9471e-03& ---&3.8423e-01& --- \\
64  &8.4301e-04&2.08&3.6239e-01&1.00&9.3626e-04&2.08&1.9102e-01&1.01 \\
128 &2.0763e-04&2.02&1.8117e-01&1.00&2.3080e-04&2.02&9.5376e-02&1.00 \\
256 &5.1709e-05&2.01&9.0582e-02&1.00&5.7493e-05&2.01&4.7671e-02&1.00 \\
512 &1.2915e-05&2.00&4.5291e-02&1.00&1.4360e-05&2.00&2.3834e-02&1.00 \\
\hline
\end{tabular}
\caption{Errors for the convergence test for \eqref{eq:solx}, with
\eqref{eq:g}, over the time interval $[0,1]$.
We also display the experimental orders of convergence (EOC).}
\label{tab:Qhxy}
\end{table}%

We begin with a convergence experiment in order to confirm our theoretical
results from Theorem~\ref{thm:main}. To this end, on recalling
\eqref{eq:xt2}, we construct the right-hand side 
\[
f_{cd} = |x_\rho|^2 x_t + y_{\rho\rho} - F_{cd}(x_\rho,y,y_\rho) y,
\]
in such a way, that the exact solution is given by a translated
and dilated circle parameterized by
\begin{equation} \label{eq:solx}
x(\rho,t)=\binom{t^2}{t^2} + (1 + t^3) \binom{\cos g(\rho)}{\sin g(\rho)},
\end{equation}
where
\begin{equation} \label{eq:g}
g(\rho) = 2\pi\rho + \delta \sin(2\pi\rho), \quad \delta = 0.1.
\end{equation}
Upon adding the correction term
\begin{equation*} 
\int_I \pi^h\left[f_{cd}(\cdot, t_{m}) \cdot \eta_h \right]\! \drho 
\end{equation*}
to the right hand side of \eqref{eq:fda}, we can perform a convergence
experiment for our scheme, comparing the obtained discrete solutions with
\eqref{eq:solx}. In particular, we define the errors 
\[
\errorxL = \max_{m=0,\ldots,M} \| x(\cdot,t_m) -  x^m_h\|_0,\quad
\errorxH = \max_{m=0,\ldots,M} \| x(\cdot,t_m) -  x^m_h\|_1,
\]
and similarly for $\erroryL$, $\erroryH$.
The obtained errors are displayed in Table~\ref{tab:Qhxy},
where we observe optimal convergence orders,
in line with the results proven in Theorem~\ref{thm:main}.
Here we partition the time interval $[0,T]$, with $T=1$, into uniform
time steps of size $\Delta t = h^2$, for $h = J^{-1} = 2^{-k}$, $k=5,\ldots,9$.
For the initial data we choose $x_h^0 = \hat x^0_h$ in line with the
assumptions of Theorem~\ref{thm:main}. But we remark that repeating the 
convergence experiment for $x_h^0 = \pi^h x^0_h$ yields very similar errors,
and the same observed orders of convergence.

\begin{figure}
\center
\includegraphics[angle=-90,width=0.4\textwidth]{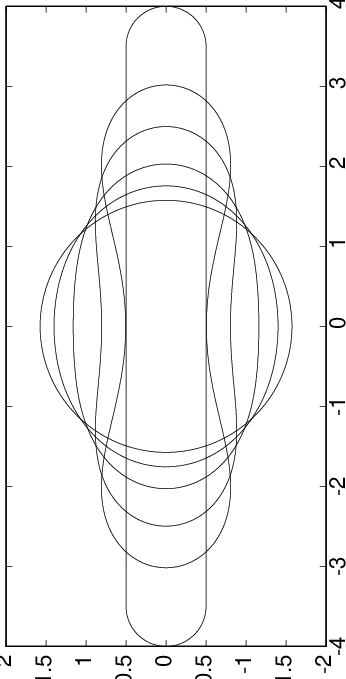}\\
\includegraphics[angle=-90,width=0.31\textwidth]{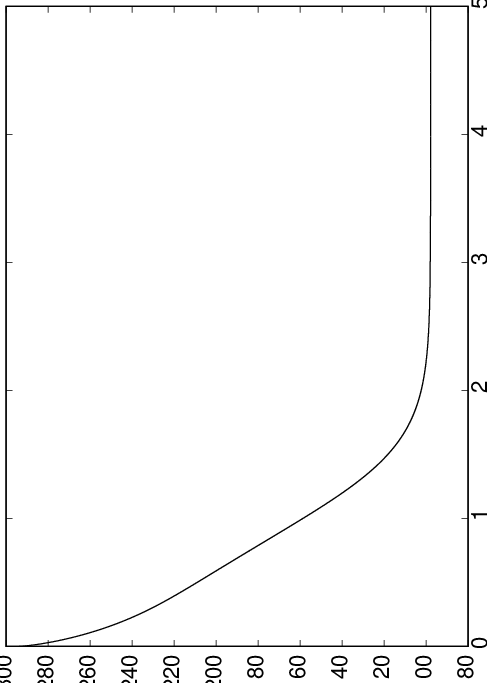}
\includegraphics[angle=-90,width=0.31\textwidth]{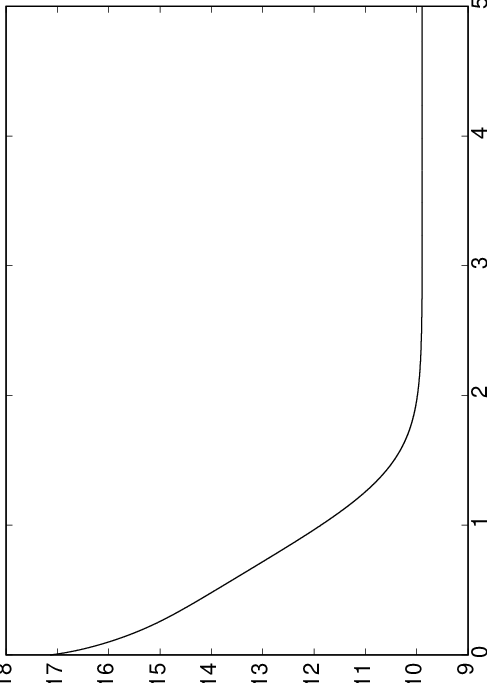}
\includegraphics[angle=-90,width=0.31\textwidth]{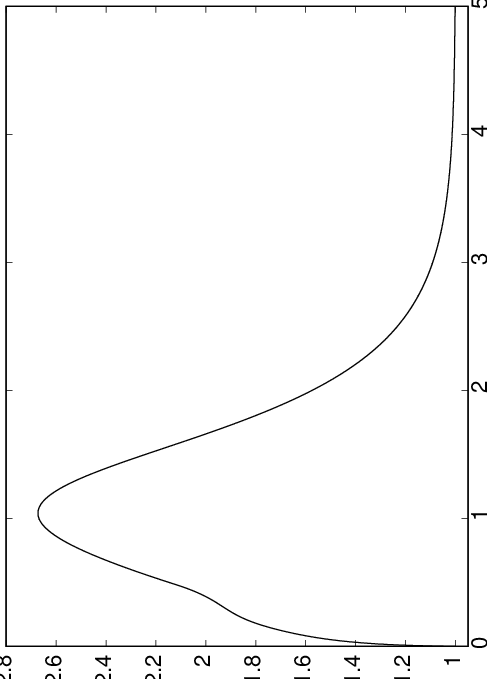}\\
{\footnotesize
\hspace{2cm} $\int_I |x^m_{h,\rho}|^2\drho$ 
\hspace{3.5cm} $|\Gamma^m|$ 
\hspace{4.2cm} $\ratio^m$} \hspace{2cm}
\caption{Curve diffusion flow for a $8\!:\!1$ tube. On top we show
$\Gamma^m$ at times $t=0,0.5,\ldots,2,T=5$. 
Below we show the evolutions of $\int_I |x^m_{h,\rho}|^2\drho$ (left),
$|\Gamma^m|$ (middle) and $\ratio^m$ (right) over time.
}
\label{fig:cigar81}
\end{figure}

The next experiment is for an elongated tube of total
dimension $8 \times 1$, and is shown in Figure~\ref{fig:cigar81}.
We can see that during the evolution the curve loses its convexity, and
eventually approaches the energy minimizing circle.
We note that the ratio \eqref{eq:ratio} 
at first increases slightly, before it eventually converges to
the optimal value of 1, which corresponds to an equidistribution
of vertices. We also observe that during the evolution shown in 
Figure~\ref{fig:cigar81}, the enclosed area was nearly preserved, 
with a relative difference of only $0.023\%$.

\begin{figure}
\center
\includegraphics[angle=-90,width=0.32\textwidth]{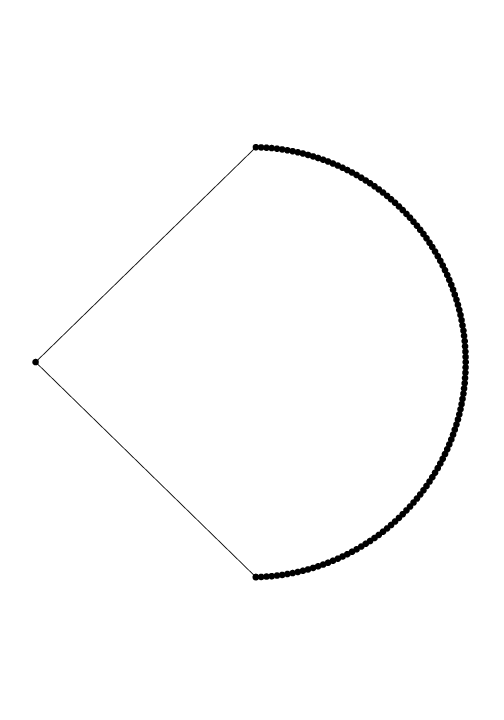}
\includegraphics[angle=-90,width=0.32\textwidth]{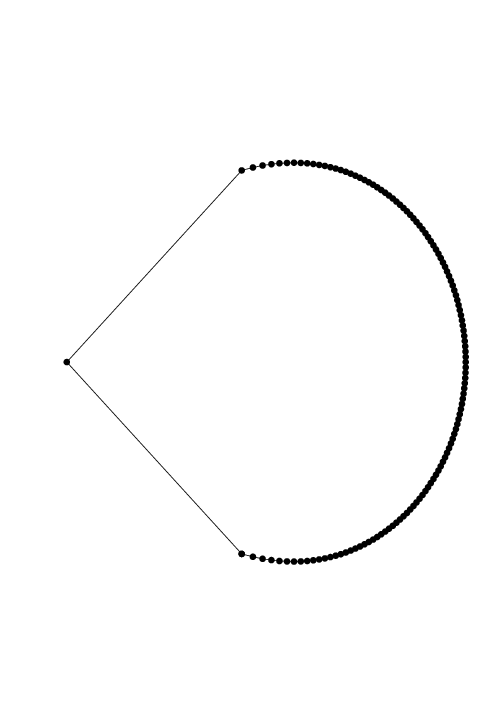}
\includegraphics[angle=-90,width=0.32\textwidth]{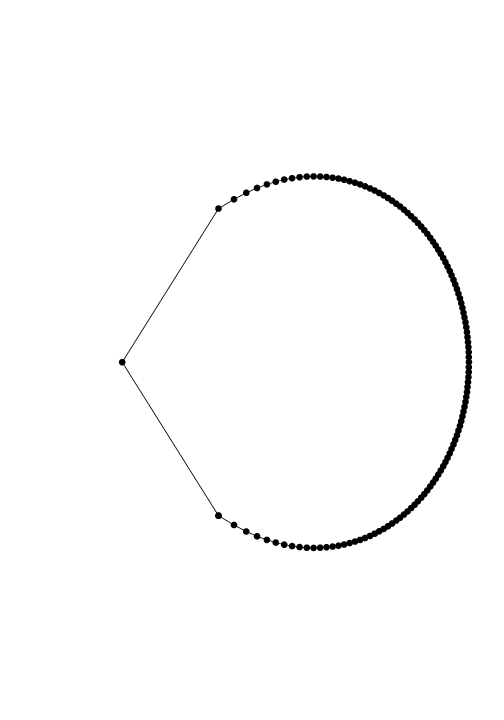}
\includegraphics[angle=-90,width=0.32\textwidth]{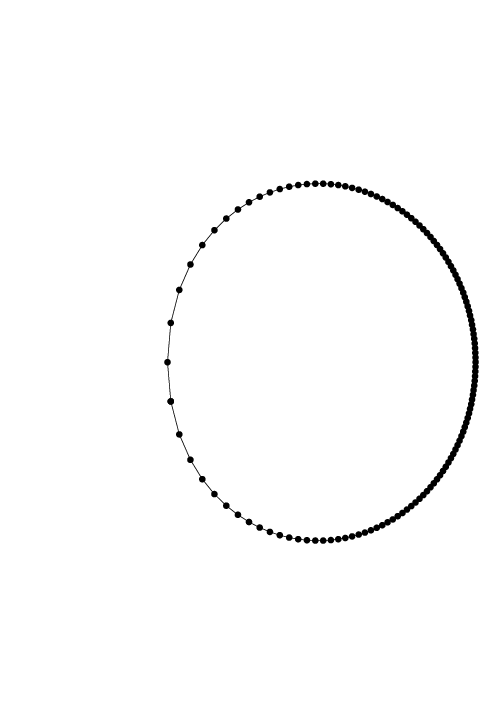}
\includegraphics[angle=-90,width=0.32\textwidth]{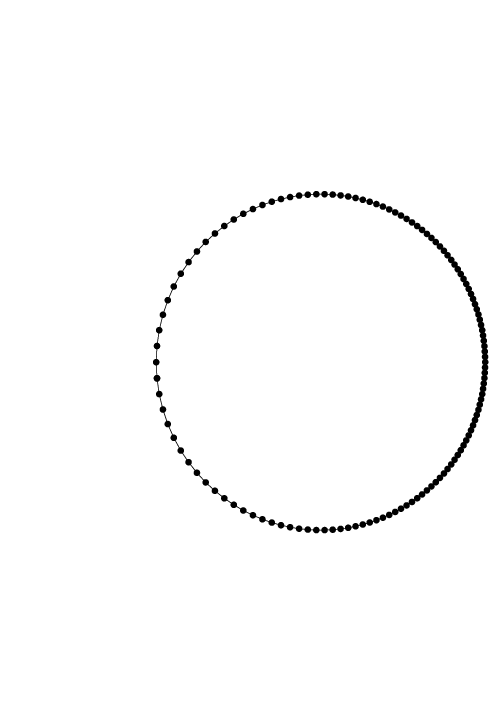}
\includegraphics[angle=-90,width=0.32\textwidth]{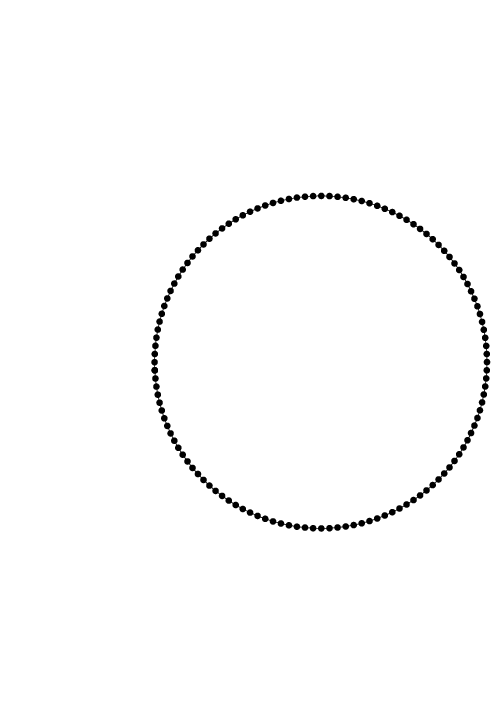}
\caption{$\Gamma^m$ at times $t=0,0.01,0.02,0.04,0.1,1$.
}
\label{fig:usp}
\end{figure}

In order to investigate the benign tangential motion further, we start our next
experiment from an initial curve that consists of a unit
semi-circle and a single additional node on the periphery of the unit circle,
compare with \cite[Fig.\ 5]{triplej}. To better highlight the movement of the
vertices, we use $J=128$ for this experiment. As can be seen from
Figure~\ref{fig:usp}, our scheme naturally moves the vertices tangentially so
that in the end an equidistributed approximation of a circle is obtained.

\begin{figure}
\center
\includegraphics[angle=-90,width=0.24\textwidth]{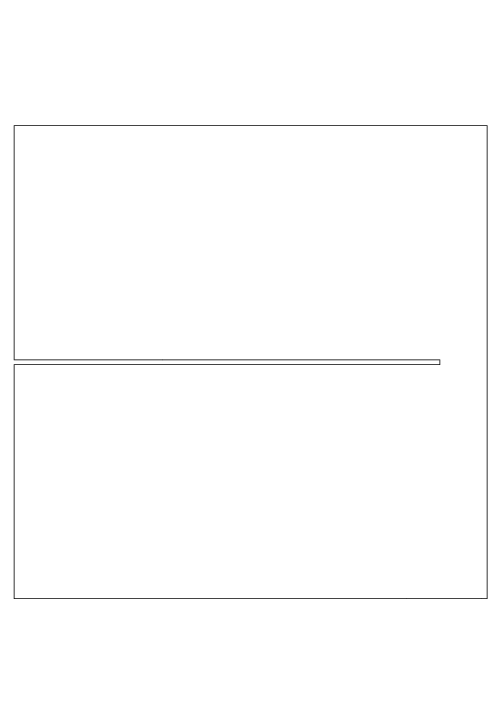}
\includegraphics[angle=-90,width=0.24\textwidth]{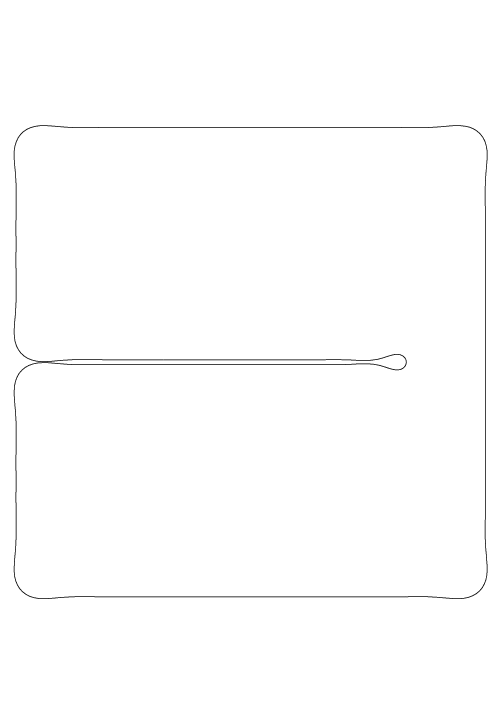}
\includegraphics[angle=-90,width=0.24\textwidth]{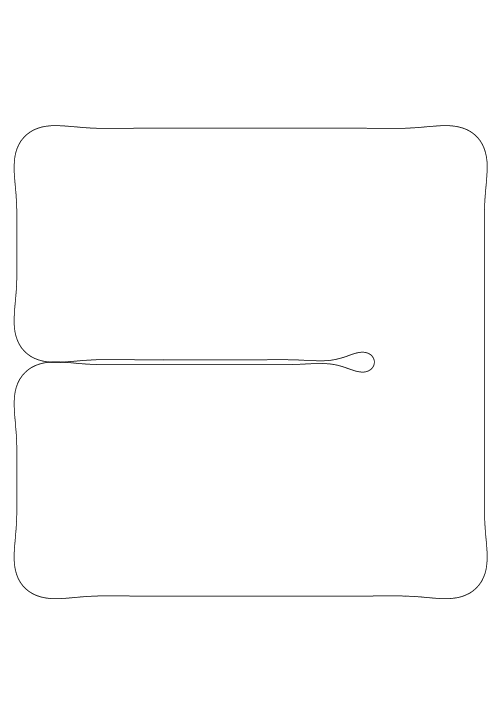}
\includegraphics[angle=-90,width=0.24\textwidth]{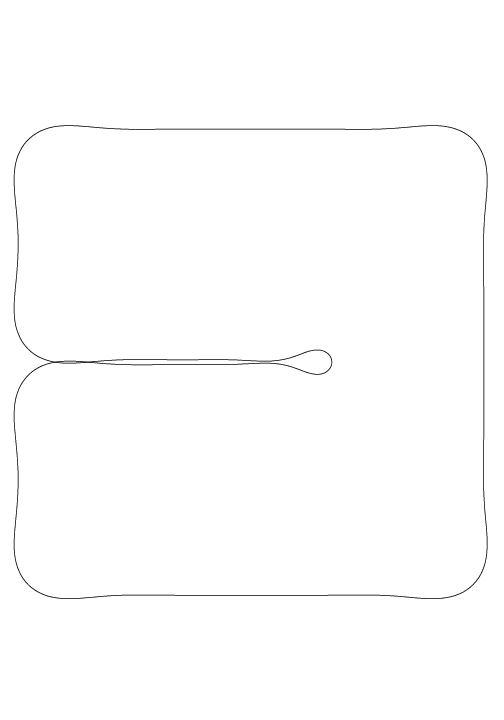}
\includegraphics[angle=-90,width=0.24\textwidth]{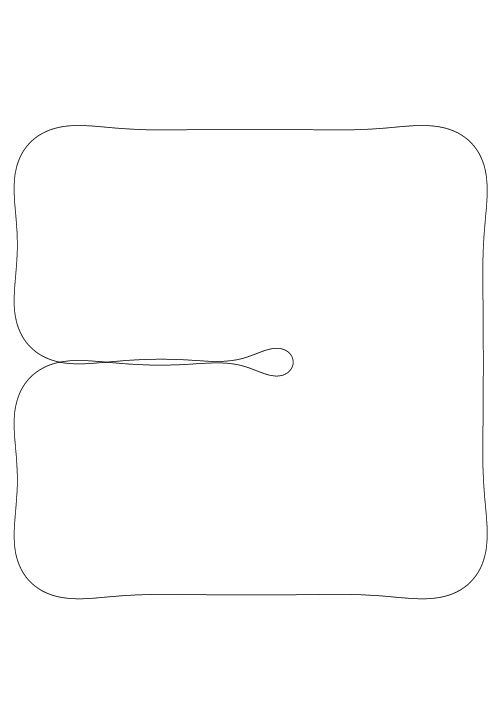}
\includegraphics[angle=-90,width=0.24\textwidth]{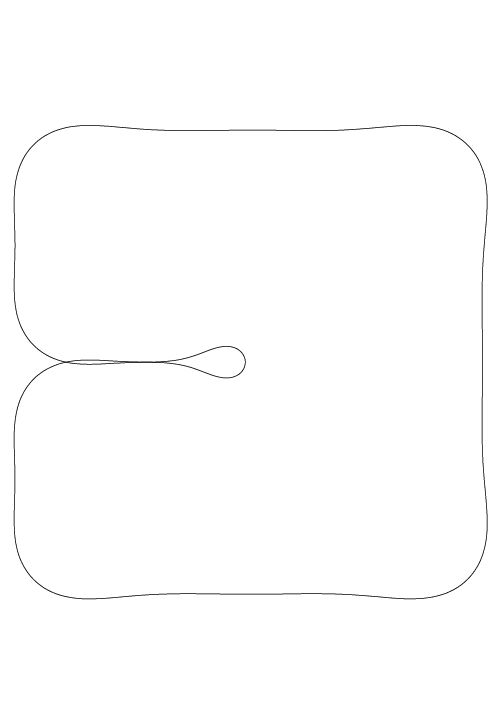}
\includegraphics[angle=-90,width=0.24\textwidth]{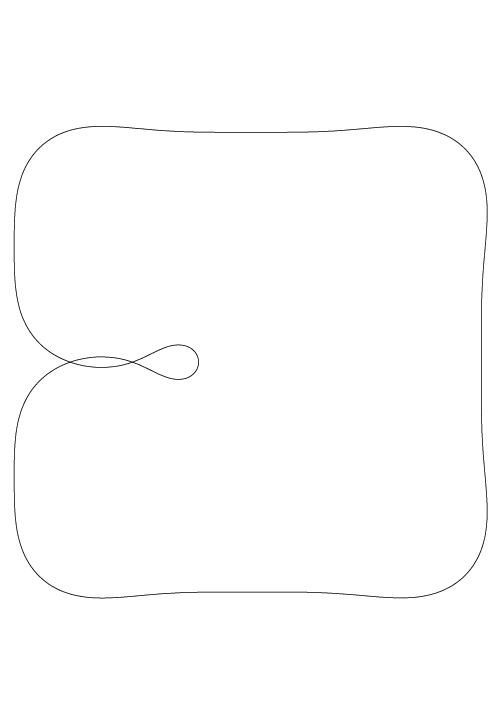}
\includegraphics[angle=-90,width=0.24\textwidth]{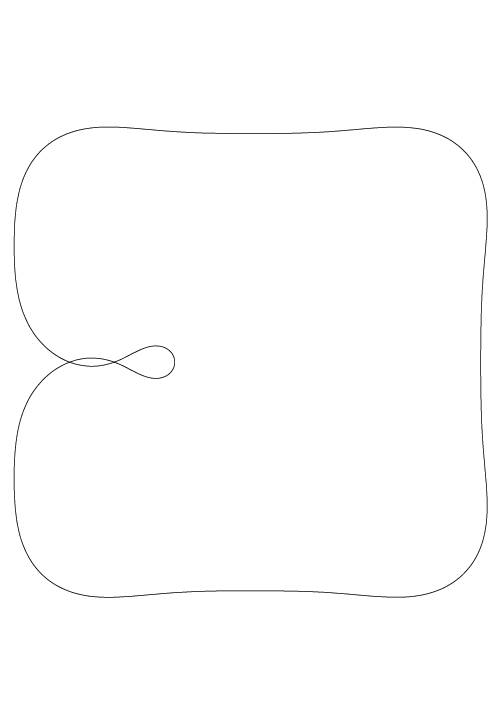}
\includegraphics[angle=-90,width=0.24\textwidth]{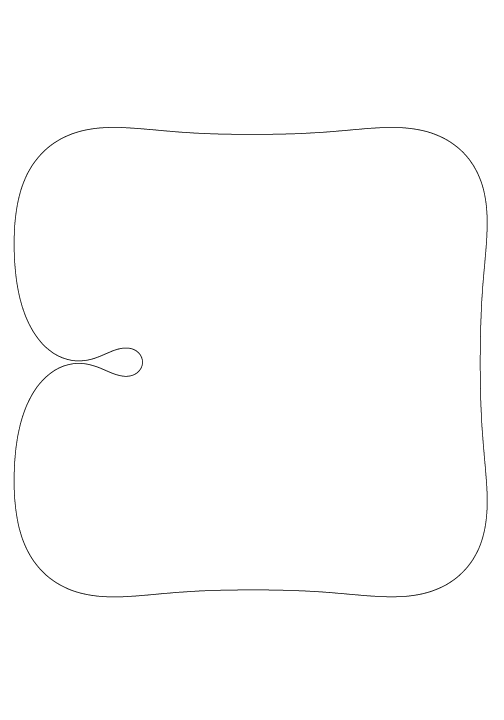}
\includegraphics[angle=-90,width=0.24\textwidth]{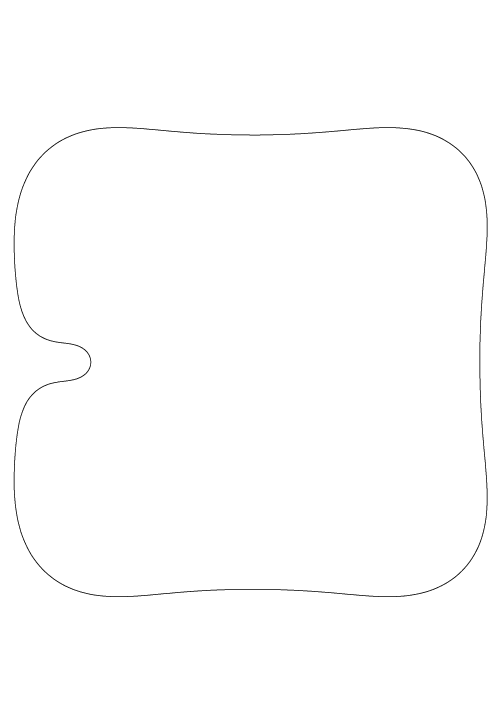}
\includegraphics[angle=-90,width=0.24\textwidth]{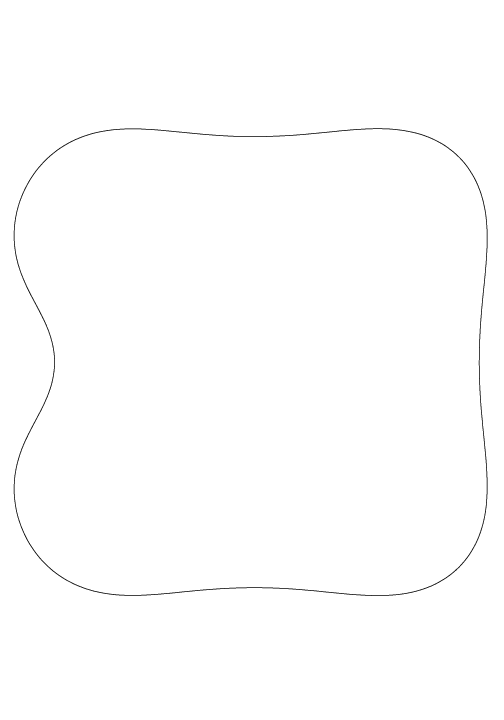}
\includegraphics[angle=-90,width=0.24\textwidth]{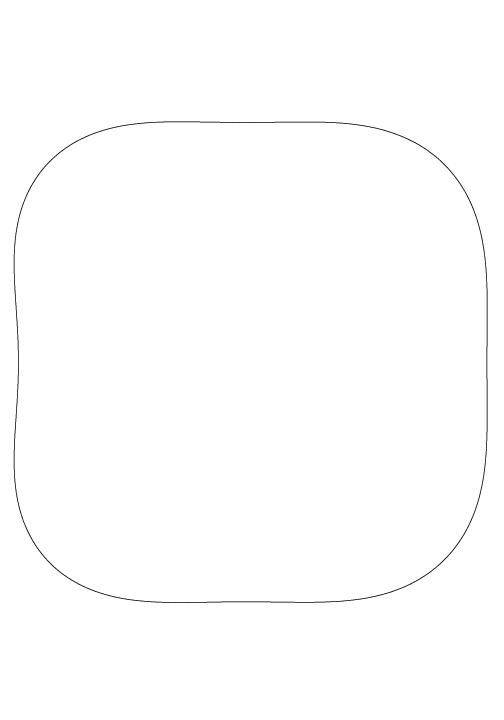}
\includegraphics[angle=-90,width=0.31\textwidth]{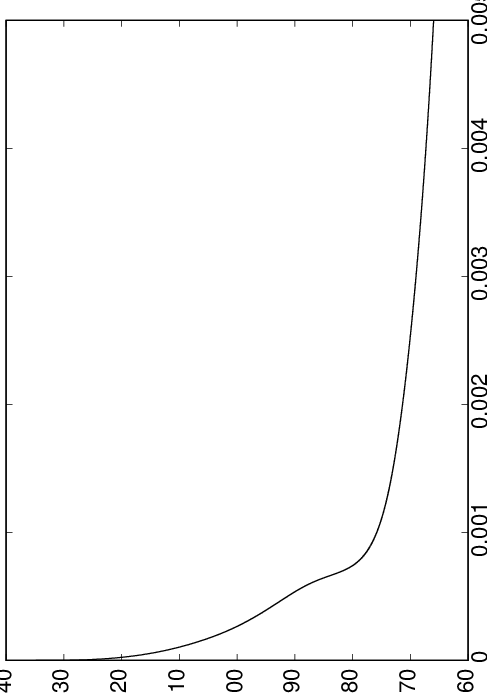}
\includegraphics[angle=-90,width=0.31\textwidth]{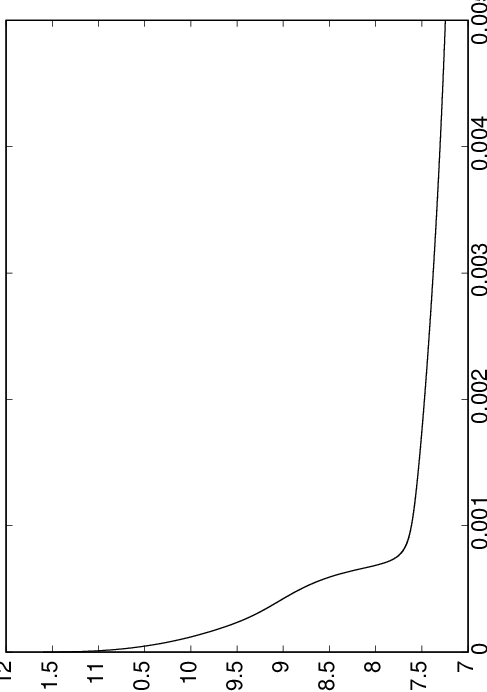}
\includegraphics[angle=-90,width=0.31\textwidth]{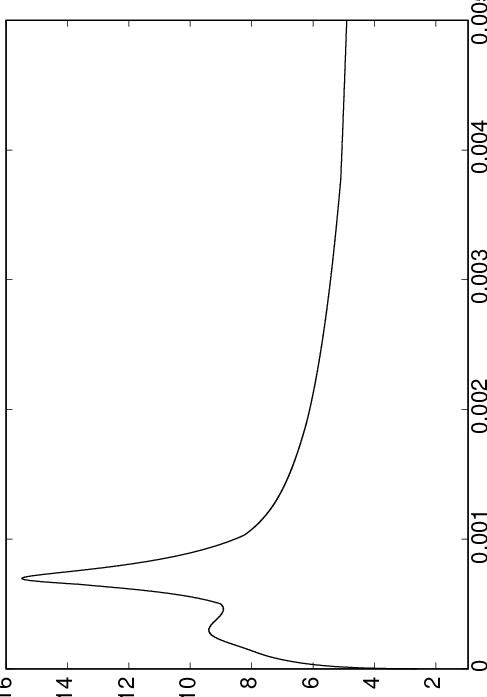}\\
{\footnotesize
\hspace{2cm} $\int_I |x^m_{h,\rho}|^2\drho$ 
\hspace{3.5cm} $|\Gamma^m|$ 
\hspace{4.2cm} $\ratio^m$} \hspace{2cm}
\caption{Curve diffusion flow for a slit domain. $\Gamma^m$ at times 
$t=0,\ 5\times10^{-6},\ 2\times10^{-5},\ 6\times10^{-5},\ 
1.2\times10^{-4},\ 2.3\times10^{-4},\ 4\times10^{-4},\ 
5\times10^{-4},\ 6\times10^{-4},\ 7\times10^{-4},\
1.1\times10^{-3},\ 
T=5\times10^{-3}$.
Below we show the evolutions of $\int_I |x^m_{h,\rho}|^2\drho$ (left),
$|\Gamma^m|$ (middle) and $\ratio^m$ (right) over time.
}
\label{fig:slit}
\end{figure}

In our final planar simulation we demonstrate that our scheme can handle 
examples with sharp corners and concavities. 
Following \cite[Fig.~16]{BanschMN05} and \cite[Fig.~9]{triplej}, we choose as 
initial data a $2\times2$ square minus a thin rectangle ($0.02\times1.8$).
Of course, the chosen initial curve does not fulfil the regularity assumptions
from Theorem~\ref{thm:main}. However, for a fixed $h$ it can be viewed as the
polygonal approximation of a smooth curve. In any case, our fully discrete 
approximation \eqref{eq:fd} can easily integrate the required evolution.
In Figure~\ref{fig:slit} we plot the results from our scheme,
using the finer time step size $\Delta t = 10^{-7}$.
The observed relative difference in area for this experiment was $0.009\%$. 
We note an excellent agreement of our results with the ones presented
in \cite[Fig.~9]{triplej}.

\begin{figure}
\center
\includegraphics[angle=-0,width=0.3\textwidth]{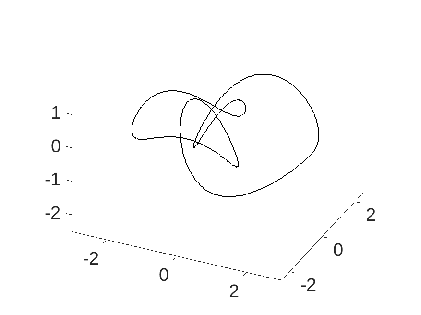}
\includegraphics[angle=-0,width=0.3\textwidth]{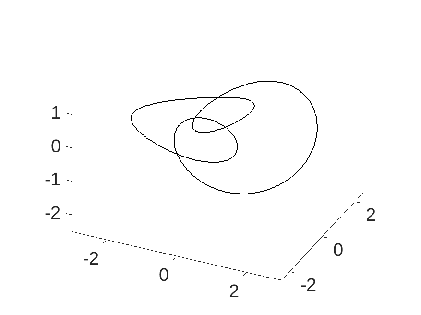}
\includegraphics[angle=-0,width=0.3\textwidth]{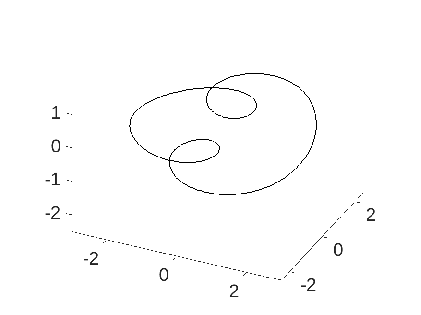}
\includegraphics[angle=-0,width=0.3\textwidth]{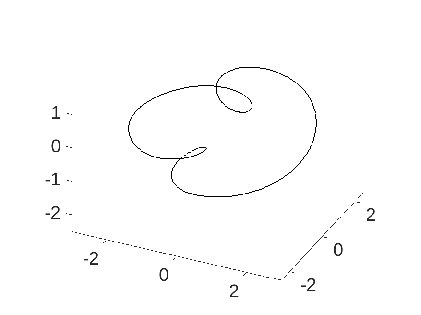}
\includegraphics[angle=-0,width=0.3\textwidth]{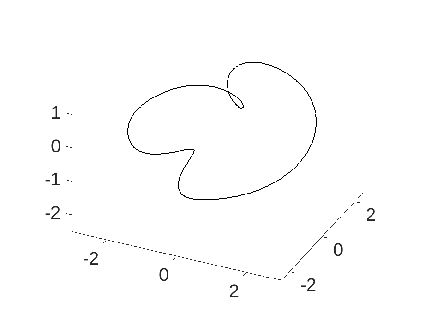}
\includegraphics[angle=-0,width=0.3\textwidth]{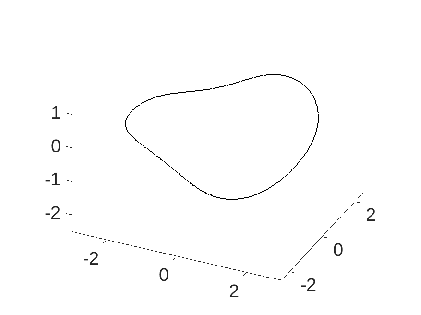}
\includegraphics[angle=-0,width=0.3\textwidth]{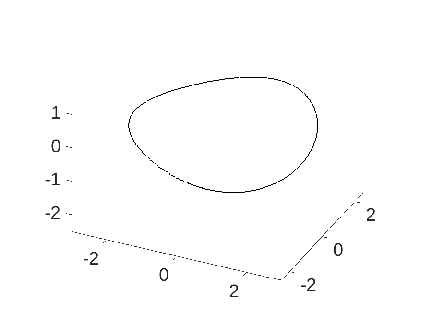}
\includegraphics[angle=-0,width=0.3\textwidth]{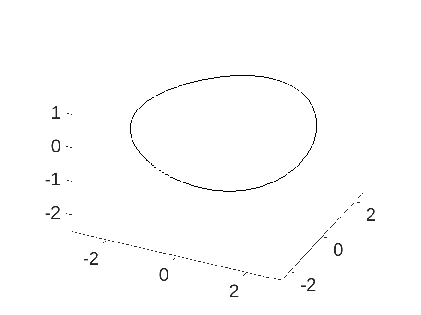}
\includegraphics[angle=-0,width=0.3\textwidth]{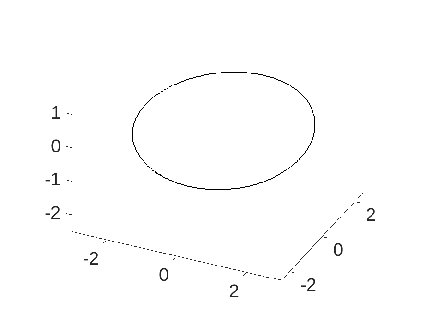} \\
\includegraphics[angle=-90,width=0.3\textwidth]{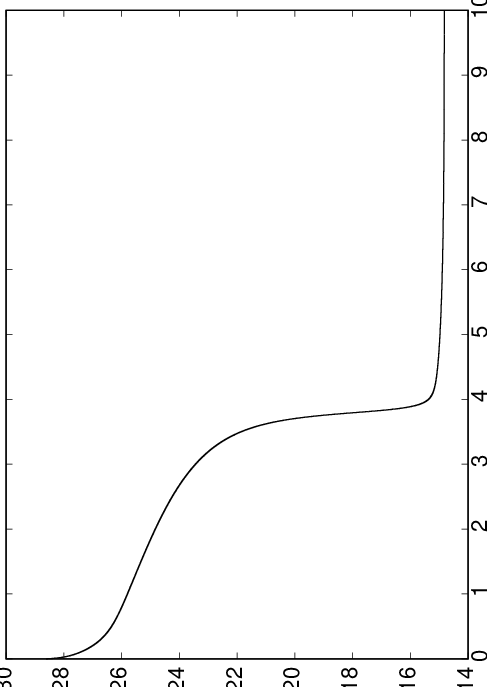}
\includegraphics[angle=-90,width=0.3\textwidth]{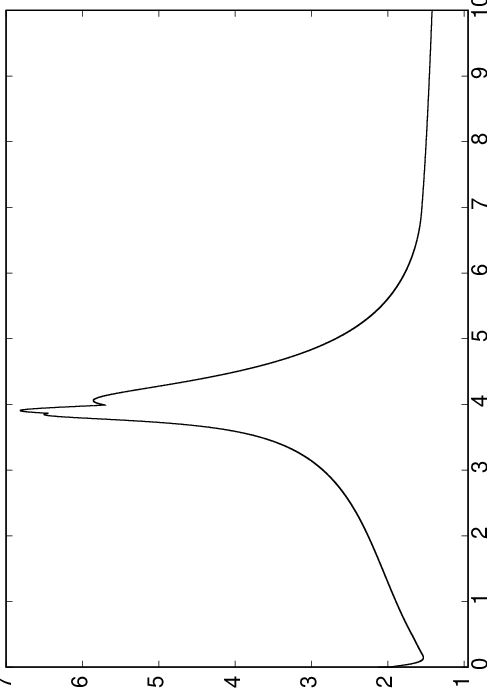}
\includegraphics[angle=-90,width=0.3\textwidth]{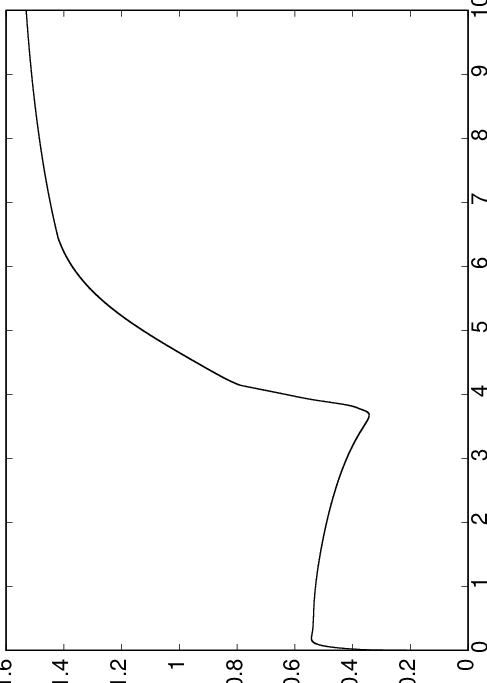}\\
{\footnotesize
\hspace{2.2cm} $|\Gamma^m|$
\hspace{4.2cm} $\ratio^m$
\hspace{4.2cm} $1/K^m_\infty$} \hspace{2cm}
\caption{Curve diffusion for the two interlocked rings \eqref{eq:irings}. 
We show $\Gamma^m$ at times $t= 0, 1, 3, 3.5, 3.7, 4, 5, 6, T=10$.
Below are plots of $|\Gamma^m|$ (left), $\ratio^m$ (middle) 
and $1/K^m_\infty$ (right) over time.}
\label{fig:irings}
\end{figure}

We also present a numerical experiment for $d=3$. To this end, we consider
the two interlocked rings in $\bR^3$ from \cite[Fig.\ 5]{fincodim}. 
In particular, the initial curve is given by
\begin{equation}
x_0(\rho) = \tfrac18 
\begin{pmatrix}
10 (\cos(2\pi\rho)+\cos(6\pi\rho))+\cos(4\pi\rho)
+\cos(8\pi\rho) \\
6\sin(2\pi\rho)+10\sin(6\pi\rho) \\
4\sin(6\pi\rho)\sin(5\pi\rho)+4\sin(8\pi\rho)-2\sin(12\pi\rho)
\end{pmatrix}
\quad \rho \in I.
\label{eq:irings}
\end{equation}
The evolution of this curve under curve diffusion can
be seen in Figure~\ref{fig:irings}. 
For this experiment we also
include a plot of the inverse of
the maximal magnitude of the discrete curvature, \eqref{eq:kappainv}. 
This strongly indicates that during the shown evolution the curve does not
become singular, i.e.\ the curvature remains bounded. This is in stark
contrast to the corresponding evolution under curve shortening flow.
In fact, the authors recently numerically studied that particular example
in \cite[Fig.\ 5]{fincodim}, where the numerical evidence suggests that
the flow undergoes a singularity. 

\subsection{Elastic flow}

In this section we will monitor the discrete elastic energy
\begin{equation*} 
E^m = \tfrac12 \int_I |P^m_h y^m_h |^2 |x^m_{h,\rho}| \drho 
+ \lambda \int_I |x^m_{h,\rho}| \drho,
\end{equation*}
where $P^m_h = \Id - \frac{x^m_{h,\rho}}{|x^m_{h,\rho}|} \otimes
\frac{x^m_{h,\rho}}{|x^m_{h,\rho}|}$. Unless otherwise stated we choose
$\lambda=0$ for the simulations presented in this section.

\begin{figure}
\center
\includegraphics[angle=-90,width=0.26\textwidth]{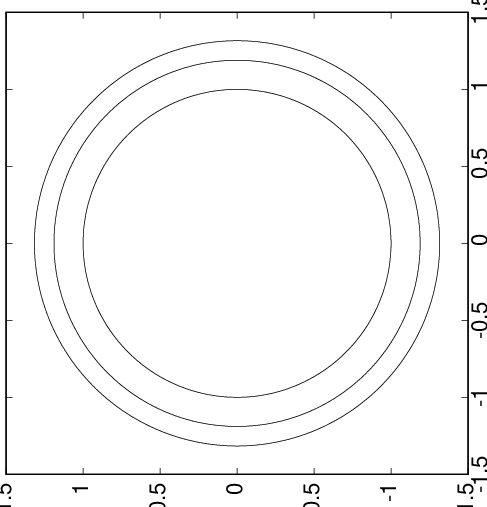}
\includegraphics[angle=-90,width=0.34\textwidth]{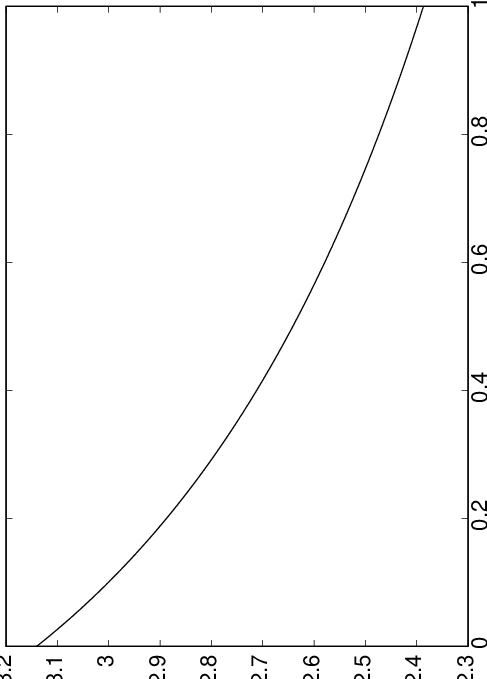}
\includegraphics[angle=-90,width=0.34\textwidth]{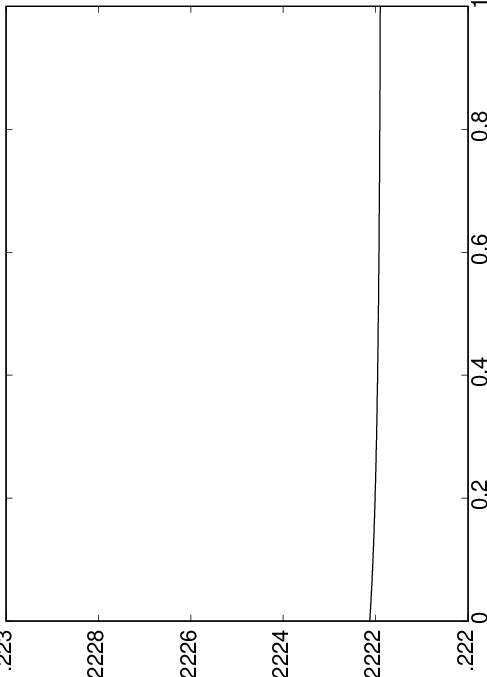} \\
{\footnotesize
\hspace{5cm} $E^m$ \hspace{5cm} $\ratio^m$}
\caption{Elastic flow with $\lambda=0$ for a circle. We show
$\Gamma^m$ at times $t=0,0.5,T=1$ (left) and plots of
$E^m$ (middle) and $\ratio^m$ (right) over time.
}
\label{fig:circle}
\end{figure}

We begin with an experiment for a known exact solution. In fact, it is easy to 
see that
\begin{equation} \label{eq:solxel}
x(\rho,t) = (1 + 2t)^{\frac14} \binom{\cos g(\rho)}{\sin g(\rho)},
\end{equation}
recall \eqref{eq:g},
satisfies \eqref{eq:Vel} with $\lambda=0$. We demonstrate the resulting
evolution by computing the numerical solution for our scheme \eqref{eq:fdel}.
As can be seen from Figure~\ref{fig:circle},
we see that the circle expands with the prescribed rate. Moreover, we observe
that our scheme induces a tangential motion that in this experiment
leads to a slightly better distribution of vertices. 
Hence the discrete functions $x_h^m$ will in general not approximate 
the particular parameterizations $x(\cdot,t_m)$, since these are strictly
radial.

Nevertheless, we can use \eqref{eq:solxel} to perform a convergence
experiment for elastic flow for our scheme \eqref{eq:fdel}. 
To this end, we construct the (tangential) right-hand side 
$f_{el} = |x_\rho|^2 x_t + y_{\rho\rho} - F_{el}(x_\rho,y,y_\rho) y$,
and then add the correction term 
$\int_I \pi^h\left[f_{el}(\cdot, t_{m}) \cdot \eta_h \right]\! \drho$
to the right hand side of \eqref{eq:fdela}. We then compare the 
obtained discrete solutions with \eqref{eq:solxel}.
The results are displayed in Table~\ref{tab:Qhxyel},
where we once again observe the expected optimal convergence rates,
similarly to Table~\ref{tab:Qhxy}.
\begin{table}
\center
\begin{tabular}{|r|c|c|c|c|c|c|c|c|}
\hline
$J$ & $\errorxL$ &EOC & $\errorxH$ &EOC & $\erroryL$ &EOC & $\erroryH$ &EOC \\ 
\hline
32  &3.5251e-03& ---&4.7691e-01& ---&4.5301e-03& ---&3.8484e-01& --- \\
64  &8.3378e-04&2.08&2.3843e-01&1.00&1.0831e-03&2.06&1.9109e-01&1.01 \\
128 &2.0533e-04&2.02&1.1921e-01&1.00&2.6766e-04&2.02&9.5385e-02&1.00 \\
256 &5.1137e-05&2.01&5.9606e-02&1.00&6.6721e-05&2.00&4.7672e-02&1.00 \\
512 &1.2772e-05&2.00&2.9803e-02&1.00&1.6668e-05&2.00&2.3834e-02&1.00 \\
\hline
\end{tabular}
\caption{
Errors for the convergence test for \eqref{eq:solxel}, with
\eqref{eq:g}, over the time interval $[0,1]$.
We also display the experimental orders of convergence (EOC).}
\label{tab:Qhxyel}
\end{table}%

\begin{figure}
\center
\includegraphics[angle=-90,width=0.3\textwidth]{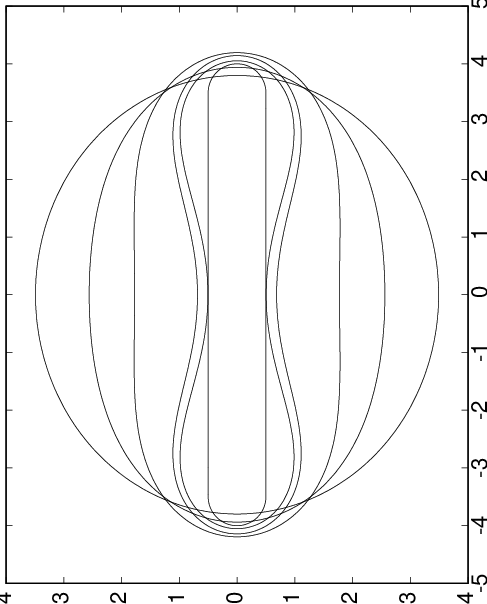}
\includegraphics[angle=-90,width=0.34\textwidth]{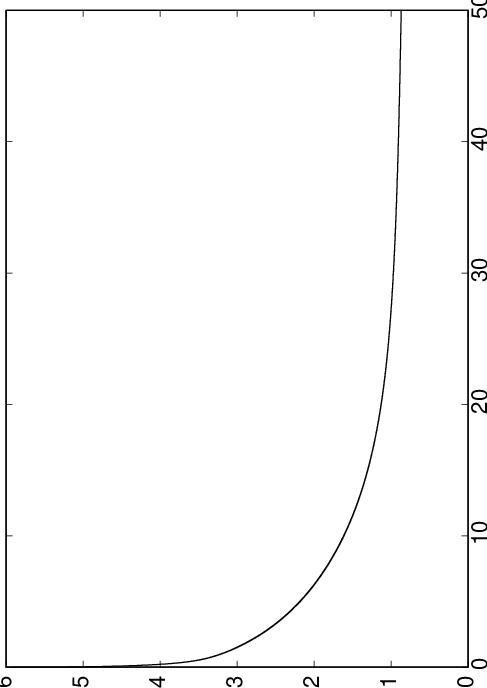}
\includegraphics[angle=-90,width=0.34\textwidth]{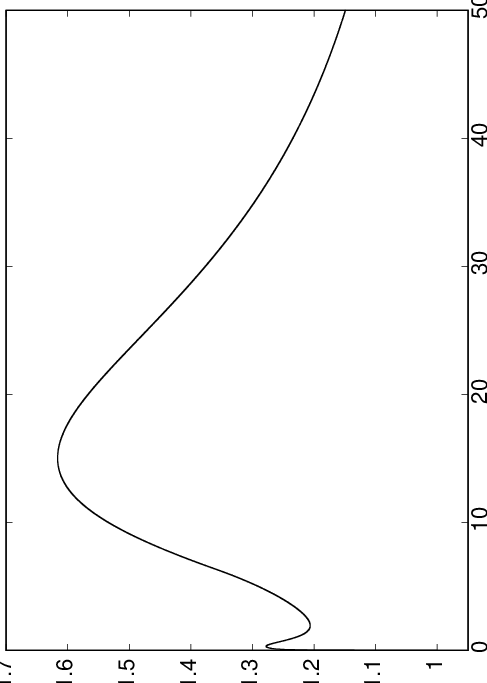} \\
{\footnotesize
\hspace{5cm} $E^m$ \hspace{5cm} $\ratio^m$}
\caption{Elastic flow with $\lambda=0$ for a $8\!:\!1$ tube. 
We show $\Gamma^m$ at times $t=0,1,2,10,20,T=50$ (left) and plots of
$E^m$ (middle) and $\ratio^m$ (right) over time.
}
\label{fig:efcigar81}
\end{figure}

We now consider an experiment for an elongated tube of total
dimension $8 \times 1$, as in Figure~\ref{fig:cigar81}.
We can see from Figure~\ref{fig:efcigar81} that during the elastic flow the
curve loses its convexity, while it monotonically decreases its discrete
energy. In fact, the curve evolves towards an expanding circle.
Moreover, the ratio \eqref{eq:ratio} at first increases slightly, before it 
eventually converges to the optimal value of 1. 
This indicates that also our scheme \eqref{eq:fdel} 
asymptotically achieves equidistribution in practice.

\begin{figure}
\center
\includegraphics[angle=-90,width=0.35\textwidth]{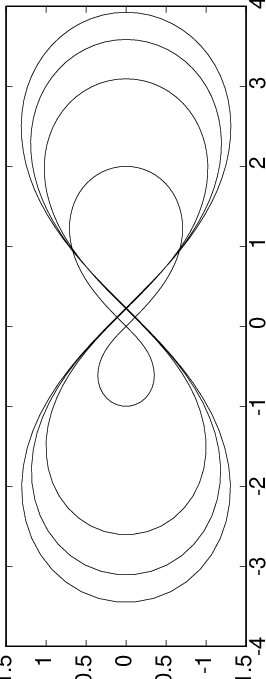}
\includegraphics[angle=-90,width=0.3\textwidth]{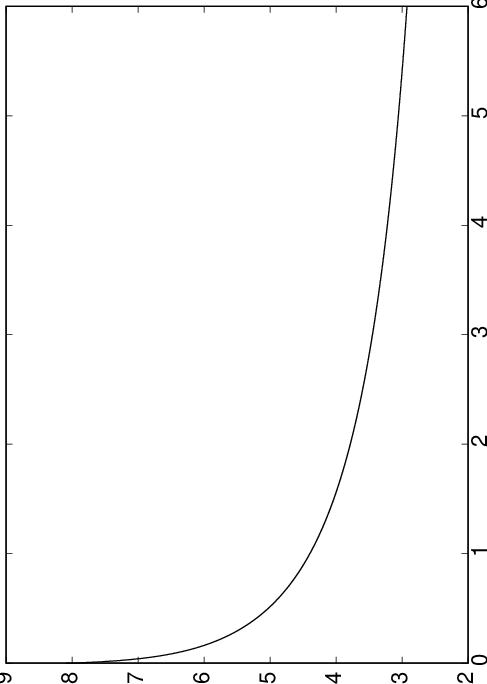}
\includegraphics[angle=-90,width=0.3\textwidth]{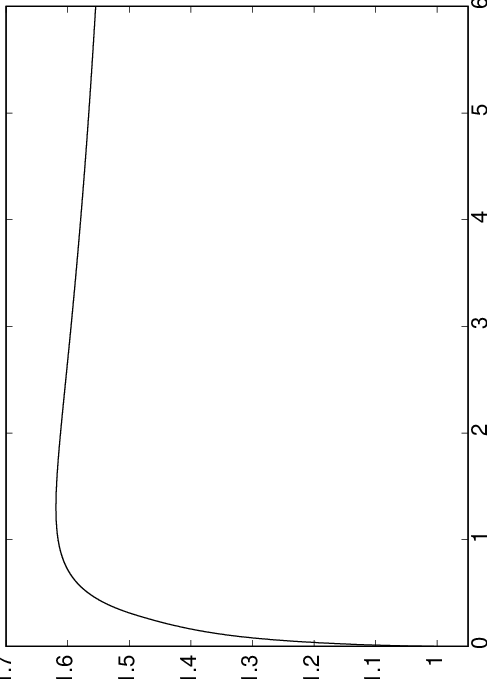} \\
{\footnotesize
\hspace{5.7cm} $E^m$ \hspace{4.5cm} $\ratio^m$}
\caption{Elastic flow with $\lambda=0$ for a $2\!:\!1$ lemniscate. 
We show $\Gamma^m$ at times $t=0,2,4,T=6$ (left) and plots of
$E^m$ (middle) and $\ratio^m$ (right) over time.
}
\label{fig:lemniscate}
\end{figure}

The next experiment is for a $2:1$ lemniscate, and is a repeat of the numerical
simulations in \cite[Fig.\ 1]{pwf}. In particular, we choose $J=100$ and 
$\Delta t=10^{-4}$ as there, and start from an equidistributed 
parameterization. We show the results obtained from our scheme
\eqref{eq:fdel} in Figure~\ref{fig:lemniscate}, and once again observe the
smooth evolution with a nice distribution of mesh points throughout. This
latter aspect is a clear difference to the corresponding results shown for
the scheme from \cite{DeckelnickD09} in \cite[Fig.\ 1]{pwf}.

\begin{figure}
\center
\includegraphics[angle=-0,width=0.3\textwidth]{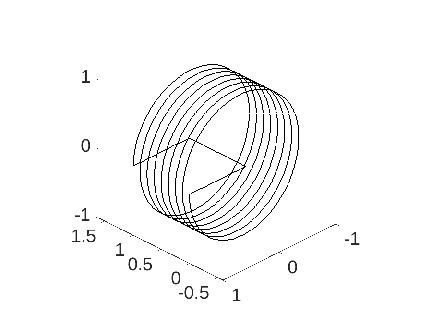}
\includegraphics[angle=-0,width=0.3\textwidth]{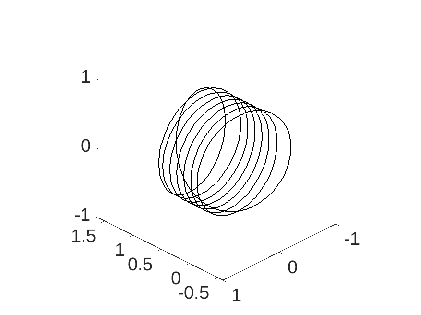}
\includegraphics[angle=-0,width=0.3\textwidth]{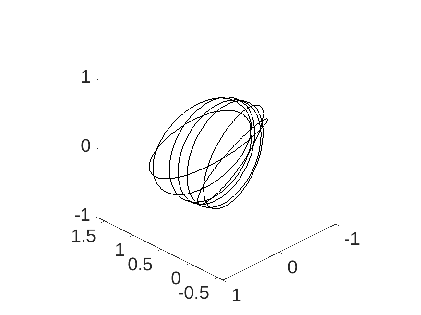}
\includegraphics[angle=-0,width=0.3\textwidth]{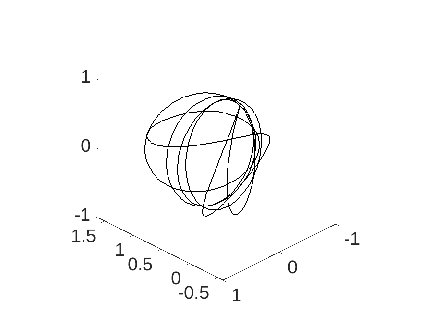}
\includegraphics[angle=-0,width=0.3\textwidth]{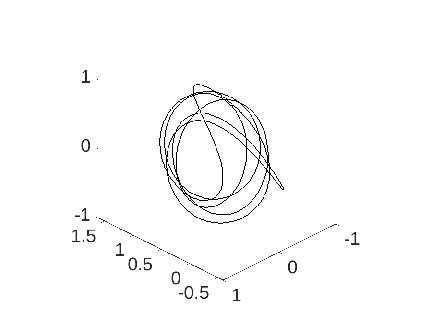}
\includegraphics[angle=-0,width=0.3\textwidth]{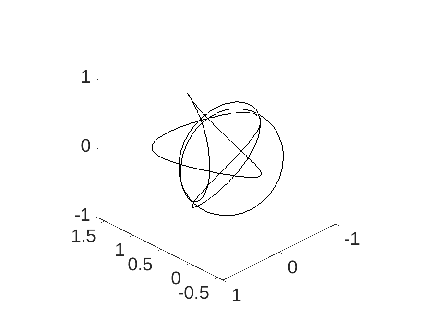}
\includegraphics[angle=-0,width=0.3\textwidth]{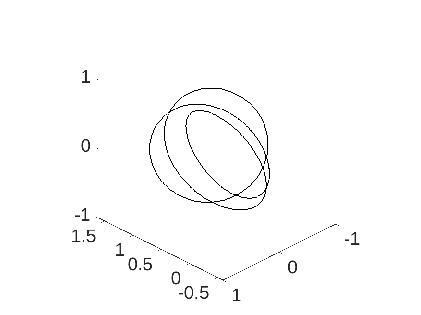}
\includegraphics[angle=-0,width=0.3\textwidth]{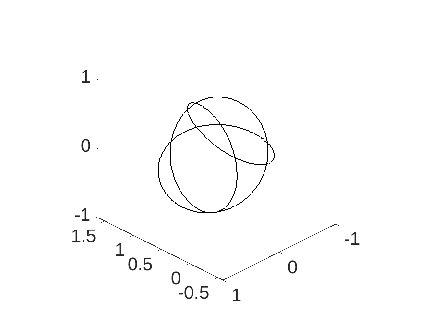}
\includegraphics[angle=-0,width=0.3\textwidth]{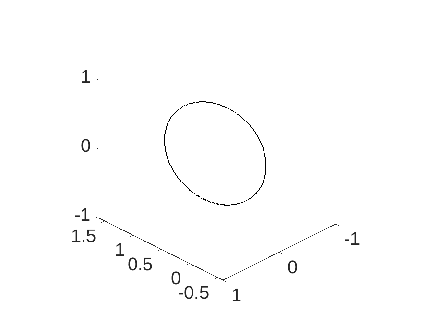}
\includegraphics[angle=-90,width=0.3\textwidth]{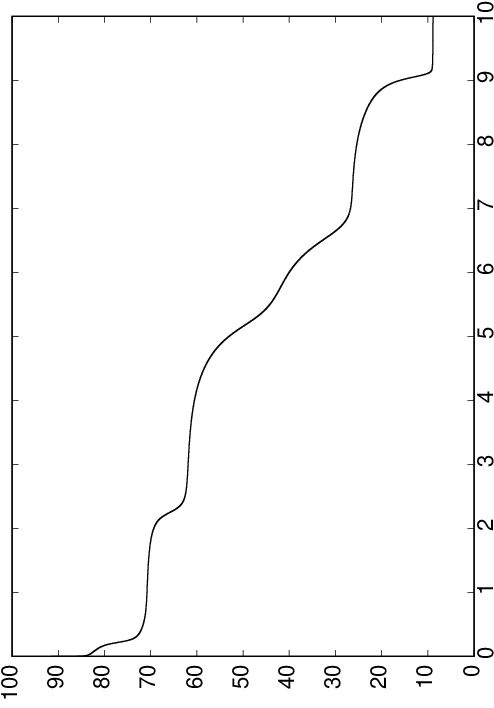}
\includegraphics[angle=-90,width=0.3\textwidth]{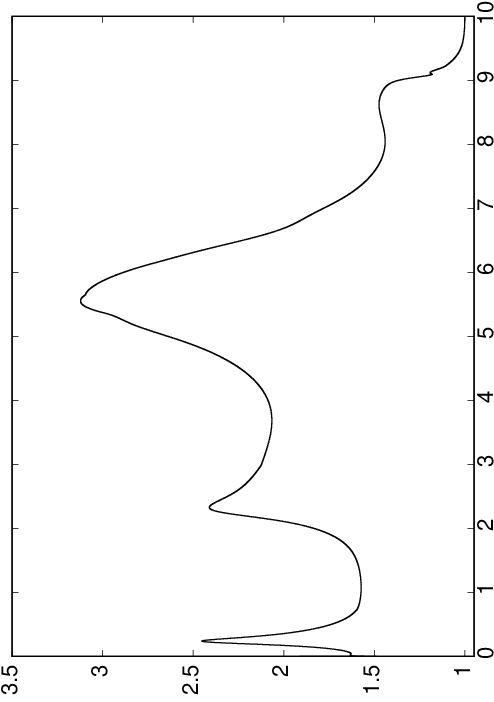} \\
{\footnotesize
\hspace{0mm} $E^m$ \hspace{4.5cm} $\ratio^m$}
\caption{Elastic flow with $\lambda=1$ for the helix \eqref{eq:helix}. 
We show $\Gamma^m$ at times $t=0,1,3,4,\ldots,8,T=10$.
Below are plots of $E^m$ (left) and $\ratio^m$ (right) over time.
}
\label{fig:el_helix}
\end{figure}

In our next experiment we consider a closed helix in $\bR^3$, 
as in \cite[Figure~2]{curves3d}. Here the open helix is defined by
\begin{equation}
 x_0(\varrho) = (\sin(16\,\pi\varrho), \cos(16\,\pi\,\varrho), \varrho)^T
\,,\quad \varrho \in [0,1]\,,
\label{eq:helix}
\end{equation}
and the initial curve is constructed from \eqref{eq:helix}
by connecting $x_0(0)$ and $x_0(1)$ with a polygon that visits the origin
and $(0,0,1)^T$. The evolution of the helix under elastic flow 
with $\lambda=1$ can be seen in Figure~\ref{fig:el_helix}. 

\begin{figure}
\center
\includegraphics[angle=-90,width=0.24\textwidth]{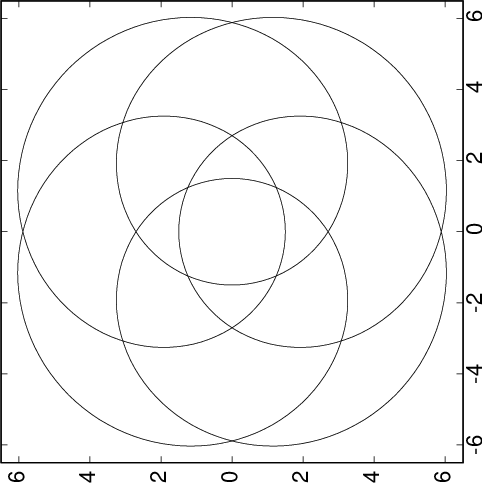}
\includegraphics[angle=-90,width=0.24\textwidth]{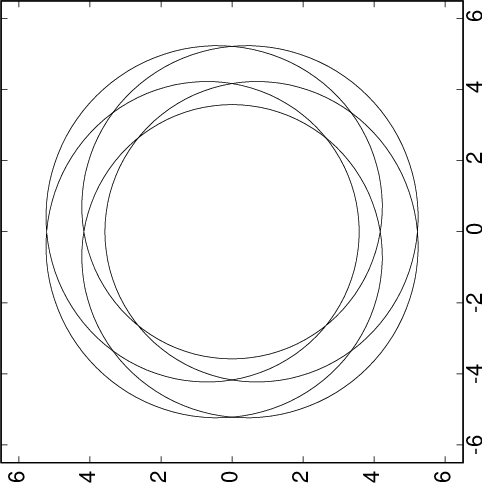}
\includegraphics[angle=-90,width=0.24\textwidth]{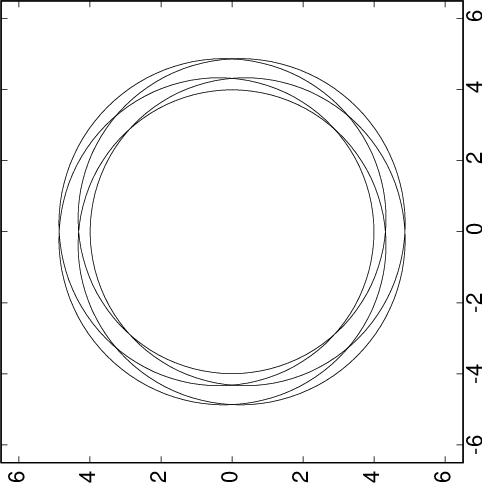}
\includegraphics[angle=-90,width=0.24\textwidth]{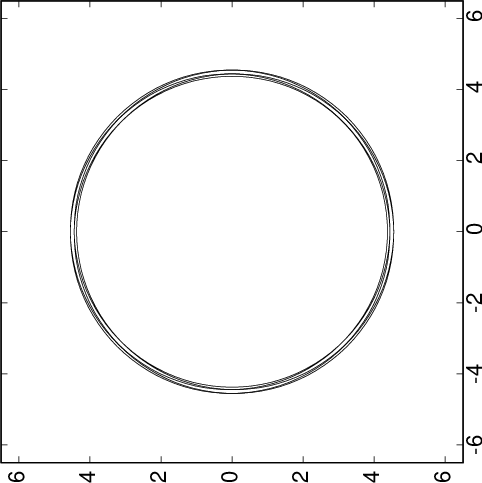}
\includegraphics[angle=-90,width=0.35\textwidth]{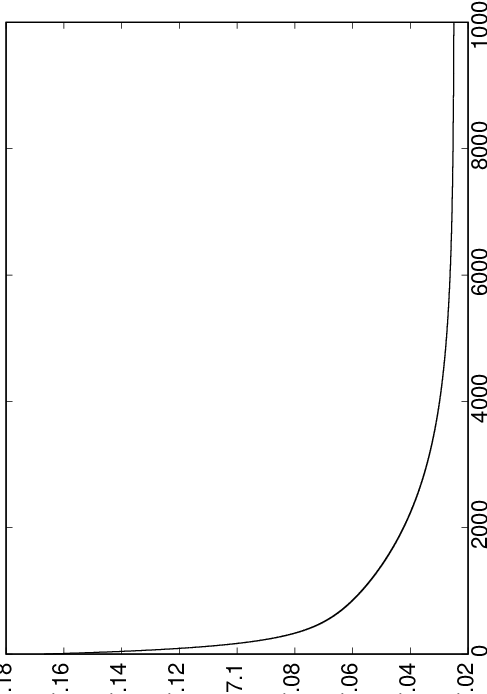}
\includegraphics[angle=-90,width=0.35\textwidth]{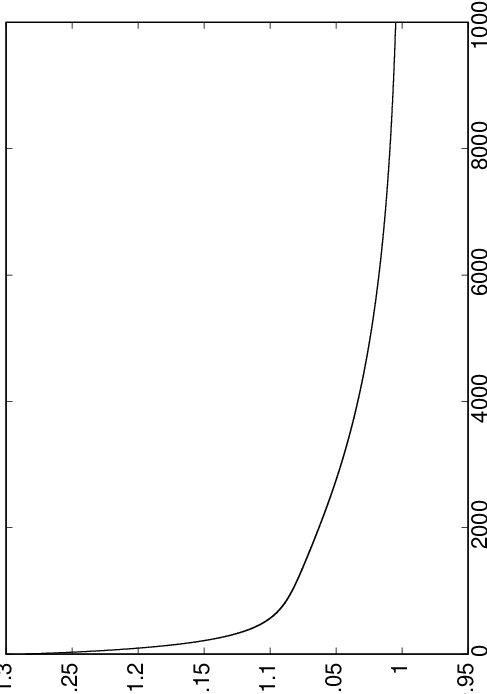} \\
{\footnotesize
\hspace{0mm} $E^m$ \hspace{5cm} $\ratio^m$}
\caption{Elastic flow with $\lambda=0.025$ for the hypocycloid
\eqref{eq:hypocycloid} with $\delta=0$. 
We show $\Gamma^m$ at times $t=0,
3000,5000,
T=10000$. 
Below are plots of $E^m$ (left) and $\ratio^m$ (right) over time.
}
\label{fig:hypocycloid}
\end{figure}

\begin{figure}
\center
\includegraphics[angle=-0,width=0.3\textwidth]{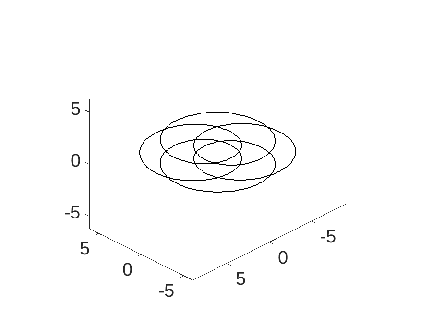}
\includegraphics[angle=-0,width=0.3\textwidth]{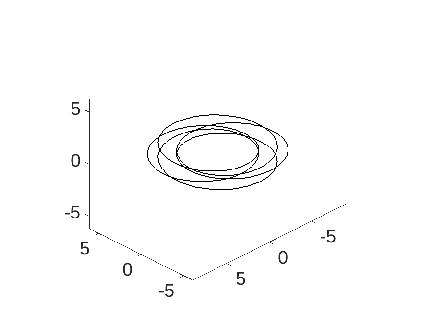}
\includegraphics[angle=-0,width=0.3\textwidth]{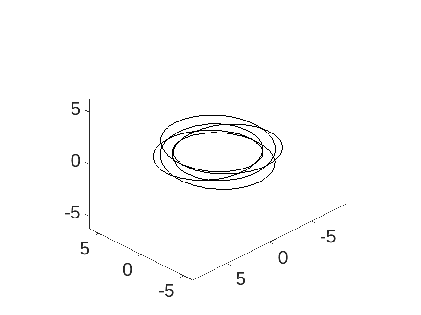}
\includegraphics[angle=-0,width=0.3\textwidth]{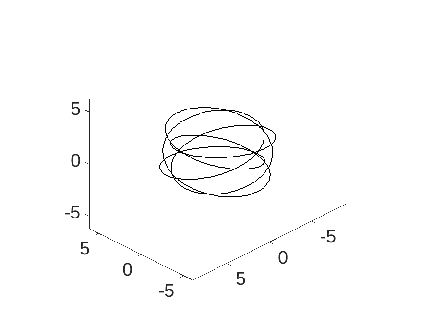}
\includegraphics[angle=-0,width=0.3\textwidth]{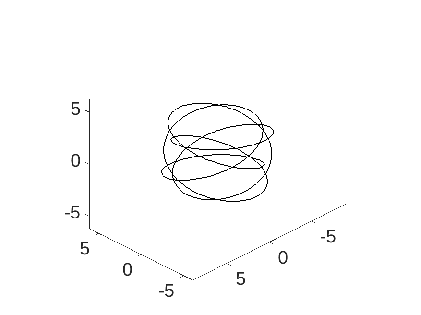}
\includegraphics[angle=-0,width=0.3\textwidth]{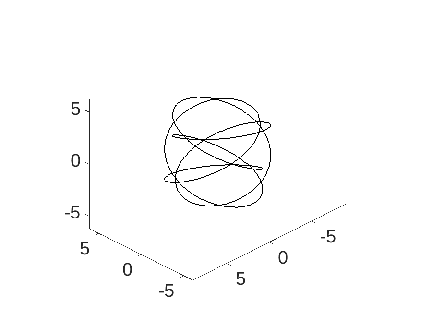}
\includegraphics[angle=-0,width=0.3\textwidth]{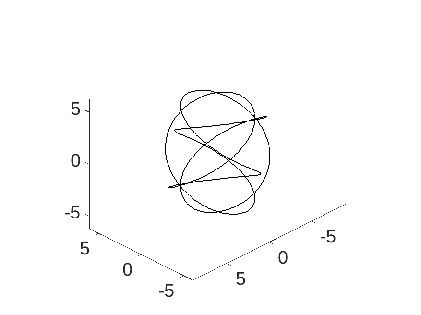}
\includegraphics[angle=-0,width=0.3\textwidth]{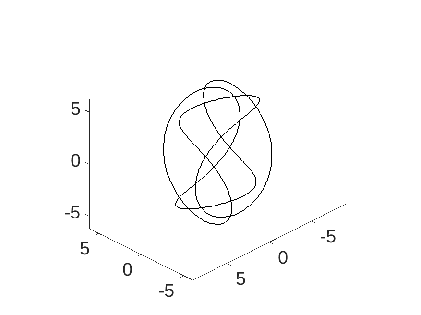}
\includegraphics[angle=-0,width=0.3\textwidth]{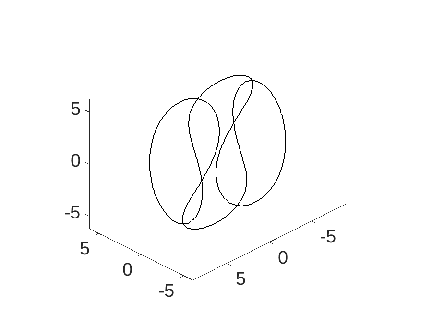}
\includegraphics[angle=-0,width=0.3\textwidth]{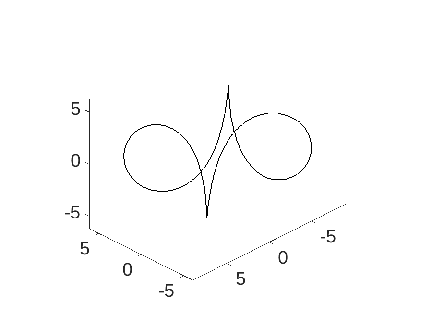}
\includegraphics[angle=-0,width=0.3\textwidth]{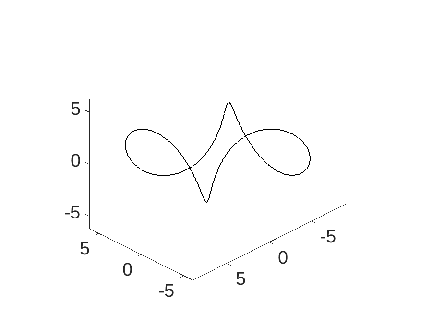}
\includegraphics[angle=-0,width=0.3\textwidth]{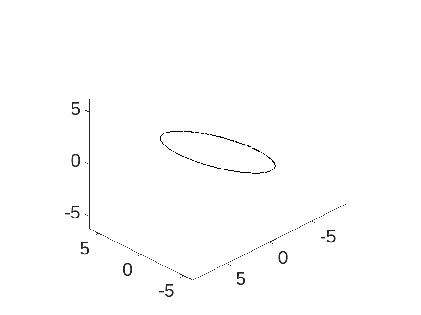}
\includegraphics[angle=-90,width=0.3\textwidth]{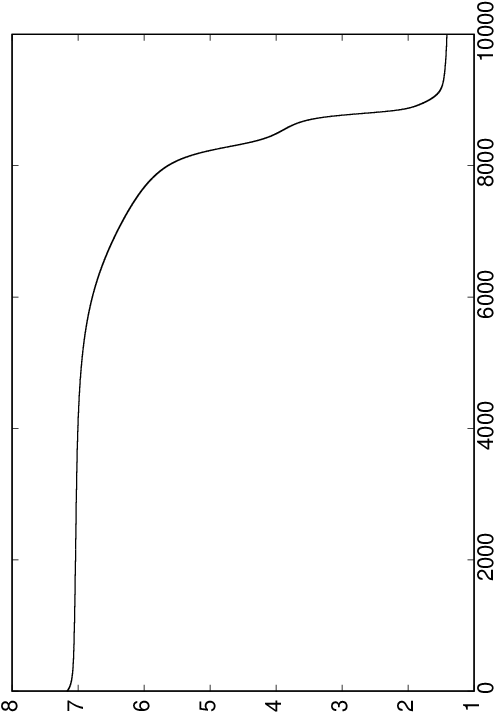}
\includegraphics[angle=-90,width=0.3\textwidth]{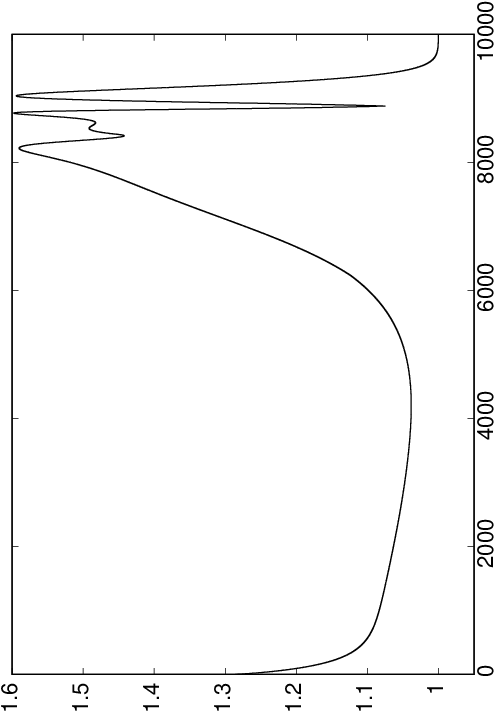} \\
{\footnotesize
\hspace{0mm} $E^m$ \hspace{4.5cm} $\ratio^m$}
\caption{Elastic flow with $\lambda=0.025$ for the hypocycloid 
\eqref{eq:hypocycloid} with $\delta=0.1$. 
We show $\Gamma^m$ at times $t=0,1700,3000,5000,5500,6000,6500,7300,8000,8500,8700,
T=10^4$.
Below are plots of $E^m$ (left) and $\ratio^m$ (right) over time.
}
\label{fig:3dhypocycloid}
\end{figure}

Our final experiments are for the hypocycloid as in 
\cite[Example~4.3]{DeckelnickD09}. In particular, we let
\begin{equation}
x_0(\rho) = \left(
-\tfrac52 \cos(2\pi\rho) + 4 \cos(10\pi\rho),
-\tfrac52 \sin(2\pi\rho) + 4 \sin(10\pi\rho),
\delta \sin(6\pi\rho)\right)
\quad \rho \in I,
\label{eq:hypocycloid}
\end{equation}
and choose $\lambda=0.025$.
For the choice $\delta=0$ the curve remains planar, evolving towards a 
five-fold covering of a circle. When we choose $\delta=0.1$,
on the other hand, the curve eventually begins to unfold and then converges to
a single circle.
For these experiments we use the larger time step size $\Delta t = 10^{-3}$. 
See Figures~\ref{fig:hypocycloid} and \ref{fig:3dhypocycloid} for the results.

\begin{appendix}
\setcounter{equation}{0} 
\renewcommand{\theequation}{\Alph{section}.\arabic{equation}}
\section{Appendix} \label{sec:AppA}

For given functions $x \in H^2(I,\bR^d), y \in H^1(I,\bR^d)$ we define 
$f: \bR^d \to \bR^d$ by
\begin{displaymath}
f(c):= F(x_\rho,y,c)y-F(x_\rho,y,0)y = 2 (x_\rho \cdot c)y +2 (x_\rho \cdot y) c - 2 (y \cdot c) x_\rho.
\end{displaymath}
Note, that the second equality holds both for $F=F_{cd}$ and $F=F_{el}$,
recall \eqref{eq:defF} and \eqref{eq:defF3}.  
Clearly, $f$ is linear and it is easily seen that  \eqref{eq:defrh}
holds if and only if $a(R_hw, \chi) = a(w,\chi)$ for all $\chi \in \Vh$, where the bilinear form 
 $a:H^1(I;\bR^d)\times H^1(I;\bR^d)
\rightarrow \bR$ is given by
\begin{displaymath}
a(v,\chi):= \int_I v_\rho \cdot \chi_\rho \drho + \int_I f(v_\rho) \cdot \chi \drho
+ \gamma \int_I v \cdot \chi \drho.
\end{displaymath}

\begin{lemma} \label{lem:Rh}
Let $x \in H^2(I,\bR^d), y \in H^1(I,\bR^d)$ with $| x_\rho | \leq C_0$, $|y | \leq C_0$ and choose $\gamma \geq 18 C_0^4+ \frac{1}{2}$. Then,
for every $w \in H^2(I;\bR^d)$ there exists a unique $R_h w \in \Vh$ such that 
\begin{displaymath}
a(R_h w, \chi) = a(w,\chi) \qquad \mbox{ for all } \chi \in \Vh
\end{displaymath}
and 
\begin{subequations} 
\begin{equation} \label{eq:Appestrha}
\| w - R_h w \|_0 + h  | w - R_h w |_1 \leq C h^2 \Vert w \Vert_2.
\end{equation}
Moreover, if $x \in C^1([0,T]; H^2(I,\bR^d)), y \in C^1([0,T];H^1(I,\bR^d))$ then $R_h w \in C^1([0,T]; \Vh)$ and we have for $t \in [0,T]$
\begin{equation} \label{eq:Appestrhb}
\| w_t(\cdot,t) -  (R_h w)_t(\cdot,t)  \|_0 +  h  | w_t(\cdot,t) - (R_h w)_t(\cdot, t) |_1 \leq C h^2 \bigl( \Vert w(\cdot,t) \Vert_2 + \Vert w_t(\cdot,t) \Vert_2 \bigr).
\end{equation}
\end{subequations}
\end{lemma}
\begin{proof}
Clearly, $| f(c) | \leq 6 C_0^2 |c|$ so that 
\begin{equation*} 
a(v,v) \geq | v |_1^2 + \gamma \| v \|_0^2 - 6 C_0^2 | v |_1 \| v \|_0 \geq \tfrac{1}{2} \Vert v \Vert_1^2,
\end{equation*}
and hence $a$ is elliptic on $H^1(I;\bR^d)$. 
The bound on $\| w - R_h w  \|_1$ follows upon inserting $\chi= \pi^h w -R_h w$ into the orthogonality relation $a(w-R_h w,\chi)=0$  and using standard interpolation estimates, 
see e.g.\ \eqref{eq:estpih}.  
The $L^2$--bound is obtained with the help of the usual  duality argument. It is not difficult to see that the 
dual problem $a(v,z)=\int_I v \cdot (w-R_h w) \drho$ for all $v \in H^1(I;\bR^d)$ is given by the equation
\begin{displaymath}
-z_{\rho \rho} - 2 \bigl( (y \cdot z) x_\rho \bigr)_\rho - 2 \bigl( (x_\rho \cdot y) z \bigr)_\rho + 2 \bigl( (x_\rho \cdot z) y \bigr)_\rho + \gamma z = w - R_h w \quad \mbox{ in } I.
\end{displaymath}
Its solution $z$ belongs to $H^2(I;\bR^d)$ and satisfies $\Vert z \Vert_2 \leq C \Vert w - R_h w \Vert_0$. The $L^2$--bound in \eqref{eq:Appestrha} is then derived in the usual way. In the case that $x$ and $y$ depend on $t$,  the structure of $f$
and integration by parts imply that $(R_hw)_t(\cdot,t)$ is a solution of the problem 
\begin{displaymath}
a((R_h w)_t(\cdot,t), \chi) = a(w_t(\cdot,t),\chi) + b \bigl( (w-R_h w)(\cdot,t),\chi \bigr)  \qquad \mbox{ for all } \chi \in \Vh,
\end{displaymath}
where the bilinear form $b$ satisfies 
\begin{displaymath}
| b \bigl( (w-R_h w)(\cdot,t),\chi \bigr) | \leq C \Vert (w-R_h w)(\cdot,t) \Vert_0 \Vert \chi \Vert_1 \leq C h^2 \Vert w(\cdot,t) \Vert_2 \Vert \chi \Vert_1
\end{displaymath}
in view of \eqref{eq:Appestrha}. The bound \eqref{eq:Appestrhb} is now obtained in a similar way as before.
\end{proof}

\begin{lemma} \label{lem:Qh}
Let $x_0\in H^4(I,\bR^d)$ such that $0<c_0 \leq |x_{0,\rho}| \leq C_0$ in $I$, set $y_0 = \frac{x_{0,\rho\rho}}{|x_{0,\rho}|^2}$ and let $\hat x^0_h \in \Vh$ be the
solution of   \eqref{eq:defQh}. Then
\begin{equation} \label{eq:Appestqh}
\| \pi^h x_0 - \hat x^0_h \|_1  \leq C h^2
\end{equation}
and there exists $h_0>0$ such that $\frac{1}{2} c_0 \leq | \hat x^0_{h,\rho} | \leq 2 C_0$ for $0<h \leq h_0$. Furthermore, for $0<h \leq h_0$ there exists a 
unique solution $\hat y_h^0 \in \Vh$ of
\begin{equation}  \label{eq:haty0App}
\int_I \hat y_h^0 \cdot \eta | \hat x^0_{h,\rho}|^2 \drho 
+ \int_I \hat x^0_{h,\rho} \cdot \eta_\rho \drho = 0
\qquad \forall\ \eta \in \Vh
\end{equation}
and we have the error bound
\begin{equation} \label{eq:Appestyh}
 \|  y_0 - \hat y^0_h \|_0 \leq C h^2.
\end{equation} 
\end{lemma}
\begin{proof}
It follows from  \eqref{eq:mvt} that 
for all $\eta \in \Vh \subset H^1(I,\bR^d)$ it holds that
\begin{displaymath}
\int_I (\pi^h x_0)_\rho \cdot \eta_\rho \drho 
= \int_I x_{0,\rho} \cdot \eta_\rho \drho =-\int_I x_{0,\rho \rho} \cdot \eta \drho
= - \int_I y_0 \cdot \eta | x_{0,\rho} |^2 \drho.
\end{displaymath}
Taking the difference of this relation with \eqref{eq:defQh} we obtain
\begin{align*}
& \int_I (\pi^h x_0 - \hat x^0_h)_\rho \cdot \eta_\rho \drho 
+ \int_I (\pi^h x_0 - \hat x^0_h) \cdot \eta \drho 
\\ & \quad
= \int_I (\pi^h y_0 - y_0) \cdot \eta | (\pi^h x_0)_\rho |^2 \drho
+ \int_I y_0 \cdot \eta \bigl[ | (\pi^h x_0)_\rho |^2 - | x_{0,\rho} |^2 \bigr] \drho 
\\ & \quad 
= \int_I (\pi^h y_0 - y_0) \cdot \eta | (\pi^h x_0)_\rho |^2 \drho \\ & \qquad
+ \int_I y_0 \cdot \eta | (\pi^h x_0 - x_0)_\rho |^2 \drho 
- 2 \int_I y_0 \cdot \eta \bigl[ x_{0,\rho} \cdot (x_0 - \pi^h x_0)_\rho \bigr] \drho.
\end{align*}
Choosing $\eta= \pi^h x_0 - \hat x^0_h$ and using integration by parts for the last integral we deduce the  bound
  \eqref{eq:Appestqh} from  \eqref{eq:estpih}. Since $c_0 \leq | x_{0,\rho} | \leq C_0$ it is in particular easily seen that there exists $h_0>0$ such that 
$\frac{1}{2} c_0 \leq | \hat x^0_{h,\rho} | \leq 2 C_0$ for $0<h \leq h_0$.\\
Next, combining \eqref{eq:defQh} and \eqref{eq:haty0App} we obtain
\begin{displaymath}
\int_I  \hat y_h^0 \cdot \eta | \hat x^0_{h,\rho} |^2 \drho =  \int_I ( \hat x^0_h- \pi^h x_0) \cdot \eta \drho + \int_I \pi^h y_0 \cdot \eta  | (\pi^h x_0)_\rho |^2  \drho
\end{displaymath}
and hence
\begin{align*}
&
\int_I ( y_0 - \hat y_h^0) \cdot \eta | \hat x^0_{h,\rho} |^2 \drho 
=
\int_I ( y_0 - \pi^h y_0) \cdot \eta | \hat x^0_{h,\rho} |^2 \drho \\
& \quad 
+ \int_I ( \pi^h x_0 - \hat x^0_h) \cdot \eta \drho 
+ \int_I \pi^h y_0 \cdot \eta \bigl[ | \hat x^0_{h,\rho} |^2 - | (\pi^h x_0)_\rho |^2 \bigr] \drho  .
\end{align*}
If we now choose $\eta = \pi^h y_0 - \hat y_h^0$ and use again the bound on $\pi^h x_0 - \hat x^0_h$ we 
deduce \eqref{eq:Appestyh} in a straightforward manner.
\end{proof}

\end{appendix}

\def\soft#1{\leavevmode\setbox0=\hbox{h}\dimen7=\ht0\advance \dimen7
  by-1ex\relax\if t#1\relax\rlap{\raise.6\dimen7
  \hbox{\kern.3ex\char'47}}#1\relax\else\if T#1\relax
  \rlap{\raise.5\dimen7\hbox{\kern1.3ex\char'47}}#1\relax \else\if
  d#1\relax\rlap{\raise.5\dimen7\hbox{\kern.9ex \char'47}}#1\relax\else\if
  D#1\relax\rlap{\raise.5\dimen7 \hbox{\kern1.4ex\char'47}}#1\relax\else\if
  l#1\relax \rlap{\raise.5\dimen7\hbox{\kern.4ex\char'47}}#1\relax \else\if
  L#1\relax\rlap{\raise.5\dimen7\hbox{\kern.7ex
  \char'47}}#1\relax\else\message{accent \string\soft \space #1 not
  defined!}#1\relax\fi\fi\fi\fi\fi\fi}


\begin{thebibliography}{10}

\bibitem{BanschDGP23}
{\sc E.~B\"{a}nsch, K.~Deckelnick, H.~Garcke, and P.~Pozzi}, {\em Interfaces:
  Modeling, Analysis, Numerics}, vol.~51 of Oberwolfach Seminars,
  Birkh\"{a}user/Springer, Cham, 2023.

\bibitem{BanschMN04}
{\sc E.~B{\"a}nsch, P.~Morin, and R.~H. Nochetto}, {\em Surface diffusion of
  graphs: variational formulation, error analysis, and simulation}, SIAM J.
  Numer. Anal., 42 (2004), pp.~773--799.

\bibitem{BanschMN05}
\leavevmode\vrule height 2pt depth -1.6pt width 23pt, {\em A finite element
  method for surface diffusion: the parametric case}, J. Comput. Phys., 203
  (2005), pp.~321--343.

\bibitem{BaoL25}
{\sc W.~Bao and Y.~Li}, {\em An energy-stable parametric finite element method
  for the planar {W}illmore flow}, SIAM J. Numer. Anal., 63 (2025),
  pp.~103--121.

\bibitem{BaoZ21}
{\sc W.~Bao and Q.~Zhao}, {\em A structure-preserving parametric finite element
  method for surface diffusion}, SIAM J. Numer. Anal., 59 (2021),
  pp.~2775--2799.

\bibitem{triplej}
{\sc J.~W. Barrett, H.~Garcke, and R.~N\"urnberg}, {\em A parametric finite
  element method for fourth order geometric evolution equations}, J. Comput.
  Phys., 222 (2007), pp.~441--462.

\bibitem{willmore}
\leavevmode\vrule height 2pt depth -1.6pt width 23pt, {\em Parametric
  approximation of {W}illmore flow and related geometric evolution equations},
  SIAM J. Sci. Comput., 31 (2008), pp.~225--253.

\bibitem{curves3d}
\leavevmode\vrule height 2pt depth -1.6pt width 23pt, {\em Numerical
  approximation of gradient flows for closed curves in {${\mathbb R}^d$}}, IMA
  J. Numer. Anal., 30 (2010), pp.~4--60.

\bibitem{fdfi}
\leavevmode\vrule height 2pt depth -1.6pt width 23pt, {\em The approximation of
  planar curve evolutions by stable fully implicit finite element schemes that
  equidistribute}, Numer. Methods Partial Differential Equations, 27 (2011),
  pp.~1--30.

\bibitem{pwf}
\leavevmode\vrule height 2pt depth -1.6pt width 23pt, {\em Parametric
  approximation of isotropic and anisotropic elastic flow for closed and open
  curves}, Numer. Math., 120 (2012), pp.~489--542.

\bibitem{bgnreview}
\leavevmode\vrule height 2pt depth -1.6pt width 23pt, {\em Parametric finite
  element approximations of curvature driven interface evolutions}, in Handb.
  Numer. Anal., A.~Bonito and R.~H. Nochetto, eds., vol.~21, Elsevier,
  Amsterdam, 2020, pp.~275--423.

\bibitem{Bartels13a}
{\sc S.~Bartels}, {\em A simple scheme for the approximation of the elastic
  flow of inextensible curves}, IMA J. Numer. Anal., 33 (2013), pp.~1115--1125.

\bibitem{DaviG90}
{\sc F.~Davi and M.~E. Gurtin}, {\em On the motion of a phase interface by
  surface diffusion}, Z. Angew. Math. Phys., 41 (1990), pp.~782--811.

\bibitem{Davis04}
{\sc T.~A. Davis}, {\em Algorithm 832: {UMFPACK} {V}4.3---an
  unsymmetric-pattern multifrontal method}, ACM Trans. Math. Software, 30
  (2004), pp.~196--199.

\bibitem{DeckelnickD95}
{\sc K.~Deckelnick and G.~Dziuk}, {\em On the approximation of the curve
  shortening flow}, in Calculus of Variations, Applications and Computations
  (Pont-\`a-Mousson, 1994), C.~Bandle, J.~Bemelmans, M.~Chipot, J.~S.~J.
  Paulin, and I.~Shafrir, eds., vol.~326 of Pitman Res. Notes Math. Ser.,
  Longman Sci. Tech., Harlow, 1995, pp.~100--108.

\bibitem{DeckelnickD09}
\leavevmode\vrule height 2pt depth -1.6pt width 23pt, {\em Error analysis for
  the elastic flow of parametrized curves}, Math. Comp., 78 (2009),
  pp.~645--671.

\bibitem{DeckelnickDE03}
{\sc K.~Deckelnick, G.~Dziuk, and C.~M. Elliott}, {\em Error analysis of a
  semidiscrete numerical scheme for diffusion in axially symmetric surfaces},
  SIAM J. Numer. Anal., 41 (2003), pp.~2161--2179.

\bibitem{DeckelnickDE05}
\leavevmode\vrule height 2pt depth -1.6pt width 23pt, {\em Computation of
  geometric partial differential equations and mean curvature flow}, Acta
  Numer., 14 (2005), pp.~139--232.

\bibitem{fincodim}
{\sc K.~Deckelnick and R.~N\"urnberg}, {\em Discrete anisotropic curve
  shortening flow in higher codimension}, IMA J. Numer. Anal., 45 (2025),
  pp.~36--67.

\bibitem{DorflerN19}
{\sc W.~D\"orfler and R.~N\"urnberg}, {\em Discrete gradient flows for general
  curvature energies}, SIAM J. Sci. Comput., 41 (2019), pp.~2012--2036.

\bibitem{DuanL24}
{\sc B.~Duan and B.~Li}, {\em New artificial tangential motions for parametric
  finite element approximation of surface evolution}, SIAM J. Sci. Comput., 46
  (2024), pp.~587--608.

\bibitem{DziukKS02}
{\sc G.~Dziuk, E.~Kuwert, and R.~Sch{\"a}tzle}, {\em Evolution of elastic
  curves in {${\mathbb R}^n$}: {E}xistence and computation}, SIAM J. Math.
  Anal., 33 (2002), pp.~1228--1245.

\bibitem{ElliottF17}
{\sc C.~M. Elliott and H.~Fritz}, {\em On approximations of the curve
  shortening flow and of the mean curvature flow based on the {D}e{T}urck
  trick}, IMA J. Numer. Anal., 37 (2017), pp.~543--603.

\bibitem{ElliottG97a}
{\sc C.~M. Elliott and H.~Garcke}, {\em Existence results for diffusive surface
  motion laws}, Adv. Math. Sci. Appl., 7 (1997), pp.~465--488.

\bibitem{EscherMS98}
{\sc J.~Escher, U.~F. Mayer, and G.~Simonett}, {\em The surface diffusion flow
  for immersed hypersurfaces}, SIAM J. Math. Anal., 29 (1998), pp.~1419--1433.

\bibitem{GageH86}
{\sc M.~Gage and R.~S. Hamilton}, {\em The heat equation shrinking convex plane
  curves}, J. Differential Geom., 23 (1986), pp.~69--96.

\bibitem{GigaI98}
{\sc Y.~Giga and K.~Ito}, {\em On pinching of curves moved by surface
  diffusion}, Commun. Appl. Anal., 2 (1998), pp.~393--405.

\bibitem{GigaI99}
\leavevmode\vrule height 2pt depth -1.6pt width 23pt, {\em Loss of convexity of
  simple closed curves moved by surface diffusion}, in Topics in nonlinear
  analysis, vol.~35 of Progr. Nonlinear Differential Equations Appl.,
  Birkh\"{a}user, Basel, 1999, pp.~305--320.

\bibitem{GoyalPL05}
{\sc S.~Goyal, N.~C. Perkins, and C.~L. Lee}, {\em Nonlinear dynamics and loop
  formation in {K}irchhoff rods with implications to the mechanics of {DNA} and
  cables}, J. Comput. Phys., 209 (2005), pp.~371--389.

\bibitem{Grayson87}
{\sc M.~A. Grayson}, {\em The heat equation shrinks embedded plane curves to
  round points}, J. Differential Geom., 26 (1987), pp.~285--314.

\bibitem{JiangL21}
{\sc W.~Jiang and B.~Li}, {\em A perimeter-decreasing and area-conserving
  algorithm for surface diffusion flow of curves}, J. Comput. Phys., 443
  (2021), p.~110531.

\bibitem{LinS04}
{\sc C.-C. Lin and H.~R. Schwetlick}, {\em On the geometric flow of {K}irchhoff
  elastic rods}, SIAM J. Appl. Math., 65 (2004), pp.~720--736.

\bibitem{MikulaS05}
{\sc K.~Mikula and D.~{\v{S}}ev{\v{c}}ovi{\v{c}}}, {\em Tangentially stabilized
  {L}agrangian algorithm for elastic curve evolution driven by intrinsic
  {L}aplacian of curvature}, in ALGORITMY 2005, A.~Handlovicova, Z.~Kriva,
  K.~Mikula, and D.~{\v S}ev{\v c}ovi{\v c}, eds., Bratislava, 2005, Slovak
  University of Technology, pp.~32--41.

\bibitem{Mullins57}
{\sc W.~W. Mullins}, {\em Theory of thermal grooving}, J. Appl. Phys., 28
  (1957), pp.~333--339.

\bibitem{Polden96}
{\sc A.~Polden}, {\em Curves and surfaces of least total curvature and
  fourth-order flows}, PhD thesis, University T{\"u}bingen, T{\"u}bingen, 1996.

\bibitem{PozziS23}
{\sc P.~Pozzi and B.~Stinner}, {\em Convergence of a scheme for an elastic flow
  with tangential mesh movement}, ESAIM Math. Model. Numer. Anal., 57 (2023),
  pp.~445--466.

\bibitem{Alberta}
{\sc A.~Schmidt and K.~G. Siebert}, {\em Design of Adaptive Finite Element
  Software: The Finite Element Toolbox {ALBERTA}}, vol.~42 of Lecture Notes in
  Computational Science and Engineering, Springer-Verlag, Berlin, 2005.

\bibitem{TaylorC94}
{\sc J.~E. Taylor and J.~W. Cahn}, {\em Linking anisotropic sharp and diffuse
  surface motion laws via gradient flows}, J. Statist. Phys., 77 (1994),
  pp.~183--197.

\bibitem{TuO-Y08}
{\sc Z.~C. Tu and Z.~C. Ou-Yang}, {\em Elastic theory of low-dimensional
  continua and its applications in bio- and nano-structures}, J. Comput. Theor.
  Nanosci., 5 (2008), pp.~422--448.

\bibitem{Wen95}
{\sc Y.~Wen}, {\em Curve straightening flow deforms closed plane curves with
  nonzero rotation number to circles}, J. Differential Equations, 120 (1995),
  pp.~89--107.

\bibitem{Wheeler73}
{\sc M.~F. Wheeler}, {\em A priori {$L\sb{2}$} error estimates for {G}alerkin
  approximations to parabolic partial differential equations}, SIAM J. Numer.
  Anal., 10 (1973), pp.~723--759.

\end{thebibliography}
\end{document}